\documentclass{amsart}

\usepackage[
paperwidth=7in,
paperheight=10in,
margin=1in
]{geometry}

\usepackage{amsmath,amssymb,amsthm,bm,cite,graphicx,algorithm,algorithmic}
\usepackage{epsfig,subfigure}
\usepackage{enumitem}
\usepackage{tablefootnote,makecell,multirow}
\usepackage{mathtools}
\usepackage{hyperref}
\usepackage{threeparttable}
\usepackage{graphicx}
\usepackage{makecell}
\usepackage{caption}

\newtheorem{theorem}{Theorem}[section]
\newtheorem{lemma}[theorem]{Lemma}

\theoremstyle{definition}
\newtheorem{definition}[theorem]{Definition}

\theoremstyle{remark}
\newtheorem{remark}[theorem]{Remark}

\numberwithin{equation}{section}

\DeclareMathOperator{\rank}{rank}

\def\twon #1{\|#1\|_2}
\def\frobn #1{\|#1\|_{\text{F}}}
\def\abs #1{|#1|}
\def\st{\text{subject to }}

\def\bR{\mathbb{R}}
\def\bE{\mathbb{E}}

\def\bS{\mathbb{S}}

\def\m #1{\boldsymbol{#1}}
\def\cA{\mathcal{A}}
\def\cB{\mathcal{B}}

\def\cE{\mathcal{E}}

\def\cI{\mathcal{I}}

\def\cL{\mathcal{L}}
\def\cM{\mathcal{M}}
\def\cN{\mathcal{N}}
\def\cO{\mathcal{O}}
\def\cP{\mathcal{P}}

\def\cT{\mathcal{T}}

\def\equ #1{\begin{equation}#1\end{equation}}   %% ref num
\def\equs #1{\begin{equation*}#1\end{equation*}}  %% no ref num
\def\sbra #1{\left(#1\right)}
\def\mbra #1{\left[#1\right]}
\def\lbra #1{\left\{#1\right\}}

\def\tr #1{\text{tr}#1}
\def\st {\text{ subject to }}

\begin{document}

\title[Rank-Optimal Sample Complexity of RGN]{Rank-Optimal Sample Complexity of Riemannian Gauss-Newton for Low-Rank Matrix Recovery}

\author[X.~Wu]{Xunmeng Wu}
\author[J.-F.~Cai]{Jian-Feng Cai}
\author[Z.~Xu]{Zongben Xu}
\author[Z.~Yang]{Zai Yang}

\subjclass[2020]{Primary 65K10; Secondary 65F55, 90C26, 15A29, 60B20}
\keywords{Low-rank matrix sensing, optimal sample complexity, quadratic convergence, finite Neumann expansion, virtual-sequence decoupling}

\date{}

\begin{abstract}
	In this paper, we study the recovery of an unknown rank-$r$ matrix $\m{X}_\star\in\bR^{n_1\times n_2}$ from Gaussian linear measurements using the Riemannian Gauss-Newton (RGN) algorithm.
	Each iteration of RGN solves the RGN equation followed by a fixed-rank retraction.
	We prove that, with high probability, RGN with spectral initialization converges quadratically to $\m{X}_\star$ using $\mathcal{O}\sbra{r(n_1+n_2)\kappa^2}$ measurements, where $\kappa$ is the condition number of $\m{X}_\star$.
	This is the first result to simultaneously achieve rank-optimal sample complexity and quadratic convergence for Gaussian low-rank matrix sensing.
	The key to the proof is a sharp spectral-norm analysis of the iterates, which requires controlling the inverse operator arising from the RGN equation.
	Our analysis develops an exact finite Neumann expansion of this inverse operator and combines it with a virtual-sequence decoupling argument and the restricted isometry property.
	Extensive numerical experiments support our theoretical findings and demonstrate the advantages of RGN.
\end{abstract}

\maketitle

\section{Introduction}
Low-rank matrix recovery aims to recover a rank-$r$ matrix
$\m{X}_\star \in \bR^{n_1 \times n_2}$ with $r<\min\lbra{n_1,n_2}$, from linear measurements
\equs{
	\m{y}
	= \cA (\m{X}_\star),
}
where $\m{y} \in \bR^m$ and the sensing operator
$\cA : \bR^{n_1 \times n_2} \to \bR^m$ is defined by
\equ{
	\mbra{\cA (\m{X})}_i
	:= \langle \m{A}_i,\m{X} \rangle,
	\qquad i
	= 1,\ldots,m. 	\label{eq:sensing}
}
Here, $\m{A}_1,\ldots,\m{A}_m \in \bR^{n_1 \times n_2}$ are known sensing matrices.
As typically $m \ll n_1n_2$, this linear inverse problem is underdetermined in the ambient matrix space. However, a rank-$r$ matrix has only $r(n_1+n_2-r)$ degrees of freedom (DoF), making exact recovery possible under suitable conditions when the sample size $m$ is of this order. 

Low-rank matrix recovery provides a common framework for a broad range of problems in machine learning, signal processing, imaging, and statistics, including matrix sensing \cite{candes2011tight,recht2010guaranteed}, matrix completion \cite{candes2012exact,cai2010singular}, phase retrieval \cite{candes2013phaselift}, spectral compressed sensing \cite{wu2024multichannel,wu2026lifted}, robust principal component analysis \cite{candes2011robust}, and quantum state tomography \cite{hsu2024quantum}.
In this paper, we focus on the Gaussian sensing model, where the entries of the sensing matrices are independent Gaussian random variables.

Two fundamental criteria for evaluating a recovery algorithm are its sample complexity and computational efficiency for exact recovery. For low-rank matrix recovery, the information-theoretic optimal sample complexity scales as $\cO(rn)$ where $n := \max\{n_1,n_2\}$, matching the order of the DoF up to constant factors. We call a sample complexity rank-optimal if it depends linearly on $rn$, while allowing dependence on the condition number $\kappa$ of $\m{X}_\star$. Computational efficiency, on the other hand, is determined by both the per-iteration cost and the iteration complexity.

A classical approach to low-rank matrix recovery is convex relaxation, which replaces the nonconvex rank function with its convex surrogate, the nuclear norm \cite{recht2010guaranteed}. Nuclear-norm minimization achieves exact recovery with the optimal sample complexity $\cO(rn)$ \cite{candes2011tight}. However, convex methods operate over the full ambient matrix space and are typically computationally expensive.

To improve computational efficiency, numerous works explicitly exploit the low-rank structure and consider the nonconvex formulation
\equ{
	\min_{\m{X} \in \bR^{n_1 \times n_2}}
	\frac{1}{2}
	\twon{ \cA (\m{X}) - \m{y} }^2,
	\qquad
	\st
	\rank(\m{X})
	= r.
	\label{eq:main_problem}
}
Two major classes of computationally efficient algorithms for \eqref{eq:main_problem} are based on low-rank matrix factorization and
Riemannian optimization, respectively.

Factorization-based methods parameterize
$\m{X} = \m{L}\m{R}^T$, where
$\m{L} \in \bR^{n_1 \times r}$ and
$\m{R} \in \bR^{n_2 \times r}$, and optimize directly over the two low-dimensional factors \cite{chi2019nonconvex}. 
Factorized gradient descent (GD) with spectral initialization converges linearly to $\m{X}_\star$ from $\cO(r^2 n \kappa^2)$ samples, with an iteration complexity of $\cO(\kappa \log(1/\varepsilon))$ \cite{tu2016low}. Its sample complexity is suboptimal in its quadratic dependence on $r$, while its convergence can be slow for ill-conditioned matrices. With a more sophisticated initialization, GD can achieve nearly $\cO(rn)$ sample complexity \cite{tu2016low}. A recent breakthrough \cite{stoger2024non} shows that, for positive semidefinite (PSD) matrix sensing, GD with spectral initialization achieves the rank-optimal sample complexity $\cO(rn\kappa^2)$, while its iteration complexity increases to $\cO(\kappa^2\log(1/\varepsilon))$. 

To remove the $\kappa$-dependence from the iteration complexity, scaled gradient descent (ScaledGD) \cite{tong2021accelerating} introduces an adaptive preconditioner into the factorized GD updates, leading to an improved iteration complexity of $\cO(\log(1/\varepsilon))$. Building on the virtual-sequence decoupling argument of \cite{stoger2024non}, the recent work \cite{li2026scaled} further establishes the rank-optimal sample complexity $\cO (rn\kappa^2)$ for ScaledGD.
Beyond first-order methods, the Gauss-Newton matrix recovery (GNMR) method \cite{zilber2022gnmr} achieves quadratic convergence, with an iteration complexity of $\cO (\log\log(1/\varepsilon))$. However, its sample complexity $\cO (r^2n\kappa)$ remains suboptimal.

Riemannian optimization offers an alternative approach to efficiently solving \eqref{eq:main_problem}. By exploiting the fact that the set of rank-$r$ matrices forms a smooth manifold, it performs optimization on the fixed-rank manifold. A representative method is Riemannian gradient descent (RGD), which, at each iteration, projects the Euclidean gradient onto the tangent space at the current iterate and then maps the update back to the manifold through a retraction. With spectral initialization, RGD converges linearly to $\m{X}_\star$ with an iteration complexity of $\cO (\log(1/\varepsilon))$ \cite{wei2016guarantees}. Building on the virtual-sequence technique of \cite{stoger2024non}, the recent work \cite{cai2025fast} establishes the rank-optimal sample complexity $\cO(rn\kappa^2)$ for RGD.

Second-order Riemannian methods can achieve substantially faster local convergence. In particular, \cite{luo2024recursive} proposes the recursive importance sketching algorithm (RISRO), which can be interpreted as solving a Riemannian Gauss-Newton equation followed by an orthographic retraction. RISRO enjoys quadratic convergence and hence an iteration complexity of $\cO (\log\log(1/\varepsilon))$. However, its theoretical sample complexity is $\cO (r^2n\kappa)$ and thus remains suboptimal. 

The above theoretical guarantees are summarized in Table~\ref{table:comparison}.
They reveal a notable gap in the existing theory.
Rank-optimal sample complexity has been established for several first-order methods, including factorized GD, ScaledGD, and RGD, whereas existing quadratically convergent methods, such as GNMR and RISRO, still require sample complexity that scales quadratically with the rank.
This naturally leads to the following question:
\begin{center}
	\emph{Can rank-optimal sample complexity and quadratic convergence be achieved simultaneously for low-rank matrix sensing?}
\end{center}

Answering this question is nontrivial. The recent rank-optimal analyses of factorized GD, ScaledGD, and RGD are tailored to their first-order updates and linear convergence. In contrast, methods with quadratic convergence typically involve more complicated updates, making the existing rank-optimal analysis not directly applicable and the improvement of their sample complexity more challenging.

\begin{table}[t]
	\centering
	\renewcommand{\arraystretch}{1.8}
	\setlength{\tabcolsep}{15pt}
	
	\begin{threeparttable}
				\captionsetup{width=0.9\textwidth}
		\caption{Comparison of nonconvex methods for low-rank matrix sensing.\tnote{a} }
		\label{table:comparison}
		
		\normalsize
		\begin{tabular}{c|c|c}
			\hline
			Methods &
			\makecell{Sample complexity} &
			\makecell{Iteration complexity} \\
			\hline \hline
			
			Factorized GD \cite{stoger2024non}\tnote{b} 
			& $r n \kappa^2$
			& $\kappa^2 \log\frac{1}{\varepsilon}$ \\
			\hline
			
			ScaledGD \cite{li2026scaled}
			& $r n \kappa^2$
			& $\log\frac{1}{\varepsilon}$ \\
			\hline
			
			RGD \cite{cai2025fast}
			& $r n  \kappa^2$
			& $\log\frac{1}{\varepsilon}$ \\
			\hline
			
			GNMR \cite{zilber2022gnmr}
			& $r^2 n  \kappa^2$
			& $\log\log\frac{1}{\varepsilon}$ \\
			\hline
			
			RISRO \cite{luo2024recursive}
			& $r^2 n  \kappa^2$
			& $\log\log\frac{1}{\varepsilon}$  \\
			\hline
			
			\makecell{\textbf{RGN} \\ (\textbf{this work})}
			& $\boldsymbol{ rn\kappa^2 }$
			& $\boldsymbol{ \log\log\frac{1}{\varepsilon} }$ \\
			\hline
		\end{tabular}
		
		\begin{tablenotes}[flushleft]
			\footnotesize
			\item[a]
			We set $n:=\max\lbra{n_1,n_2}$ and omit the $\cO(\cdot)$ notation for simplicity. All methods use spectral initialization.
			\item[b] The guarantee for factorized GD is established for the PSD matrix sensing.
		\end{tablenotes}
		
	\end{threeparttable}
\end{table}

\subsection{Our Contribution}
In this paper, we give an affirmative answer to the question above. We study a Riemannian Gauss-Newton (RGN) method for low-rank matrix sensing. Each iteration of RGN consists of solving the RGN equation on the tangent space and applying the rank-$r$ truncated singular value decomposition (SVD) as a retraction.

We prove that, under Gaussian sensing and spectral initialization, RGN converges quadratically to $\m{X}_\star$ with high probability, provided that the sample size satisfies
\equs{
	m = \cO(r(n_1+n_2)\kappa^2),
}
with an iteration complexity of $\cO(\log\log(1/\epsilon))$. Thus, RGN achieves both rank-optimal sample complexity and quadratic convergence.

The main technical challenge is to obtain a sharp spectral-norm bound for the RGN update. The solution to the RGN equation involves the inverse of a tangent-space operator, for which the standard restricted-isometry argument naturally yields a Frobenius-to-Frobenius bound; see Lemma \ref{lem:tangent_RIP}. Directly converting this bound to the spectral norm would lead to a suboptimal sample complexity.

To overcome this difficulty, we introduce an exact finite Neumann expansion of the inverse operator; see \eqref{eq:bound1}. This expansion decomposes the inverse into a finite sum of inverse-free terms and an exact remainder. The finite terms can be sharply controlled in the spectral norm via a tailored virtual-sequence decoupling argument inspired by \cite{stoger2024non,cai2025fast}.
The remainder still contains the inverse operator, but is multiplied by a high power of a contraction operator. We therefore apply the available Frobenius-norm bound only to the remainder. By choosing the expansion order sufficiently large, the resulting geometric decay compensates for the additional $\sqrt{r}$ factor introduced by the Frobenius-norm bound.

\subsection{Related Work}
As a factorization-based method, alternating minimization \cite{jain2013low} alternately solves least-squares subproblems for the two factors. Its theoretical sample complexity has a cubic dependence on the rank, and only linear convergence is established.

Other Riemannian optimization methods for low-rank matrix recovery include Riemannian conjugate gradient (RCG) \cite{wei2016guarantees} and preconditioned RGD (PRGD) \cite{bian2024preconditioned}. Their sample complexities exhibit a suboptimal quadratic dependence on the rank, and their convergence rates remain linear.

RGN methods have also been studied in general manifold optimization, low-rank tensor estimation, and low-rank matrix completion. The classical RGN framework is presented in \cite{absil2008optimization}, with its local convergence properties further analyzed in \cite{breiding2018convergence}. For low-rank tensor estimation, \cite{luo2023low} established recovery guarantees together with quadratic convergence. However, when specialized to matrix sensing, their recovery guarantee does not yield rank-optimal sample complexity.

For low-rank matrix completion, a recent work \cite{cai2026near} showed that RGN with a multiscale residual initialization achieves exact recovery from $\mathcal{O}(\mu rn\log n\log(2\mu r\kappa))$ observed entries and eventually converges quadratically, where $\mu$ is the incoherence parameter.
Their analysis uses a leave-one-out argument to control the recovery error in the sharp norm. In contrast, we study RGN for Gaussian matrix sensing with spectral initialization. Our analysis relies on a finite Neumann expansion and a virtual-sequence decoupling argument.

Other nonconvex algorithms for low-rank matrix recovery include singular value projection (SVP) \cite{jain2010guaranteed}, normalized iterative hard thresholding (NIHT) \cite{tanner2013normalized}, and iteratively reweighted least squares (IRLS) \cite{fornasier2011low,mohan2012iterative,kummerle2018harmonic}. These methods operate on full-dimensional matrices and can therefore be computationally expensive for large-scale problems.

\subsection{Notation}
Let $\bR$ denote the real numbers.
Let $\bS^{n-1}:=\lbra{\m{x}\in\bR^n:\twon{\m{x}}=1}$ denote the unit sphere in $\bR^n$.
For a vector $\m{x}$, $\twon{\m{x}}$ denotes its $\ell_2$-norm.
For a matrix $\m{X}$, $\twon{\m{X}}$ and $\frobn{\m{X}}$ denote its spectral and Frobenius norms, respectively. For a rank-$r$ matrix $\m{X}$, its singular values are denoted by $\sigma_1(\m{X}) \ge \cdots \ge \sigma_r(\m{X})$.
For a linear operator $\cA: \bR^{n_1 \times n_2} \to \bR^{m}$, we denote its adjoint by $\cA^*$ and its operator norm by $\twon{\cA}:=\sup_{\frobn{\m{X}}=1}\twon{\cA(\m{X})}$.
The identity operator is denoted by $\cI$.
We use $(\cdot)^T$ to denote the transpose, and $\cN(0,1)$ to denote the standard Gaussian distribution. The standard inner product is defined as
$\langle \m{A},\m{B}\rangle:=\tr(\m{A}\m{B}^T)$.
Finally, $m\gtrsim n$ and $m\lesssim n$ mean $m\geq Cn$ and $m\leq Cn$, respectively, for some universal constant $C>0$.

\subsection{Organization}
The rest of the paper is organized as follows. Section \ref{sec:alg} introduces the RGN algorithm. Section \ref{sec:theory} presents our main theoretical result, Theorem \ref{thm:main}, outlines the main proof ideas, and provides its proof. Section \ref{sec:experiments} describes the implementation of RGN and presents numerical experiments to demonstrate its empirical performance. Section \ref{sec:conclusion} concludes the paper and discusses potential directions for future work. The appendix provides technical details for the proof of Theorem \ref{thm:main}.

\section{Algorithm} \label{sec:alg}

\subsection{Preliminaries}
\subsubsection{Fixed-rank matrix manifold}
\noindent

Define $\cM_r := \lbra{ \m{X}\in\bR^{n_1\times n_2}: \rank(\m{X}) = r }$.
The following geometric structure of the fixed-rank manifold is standard; see, e.g., \cite{vandereycken2013low}.
The set $\cM_r$ is a smooth embedded submanifold of dimension $r(n_1+n_2-r)$.
We equip $\cM_r$ with the Riemannian metric induced by the Frobenius inner product.
Let $\m{X} = \m{U}\m{\Sigma}\m{V}^T$ be a compact SVD, where $\m{U} \in \bR^{n_1 \times r}$ and $\m{V} \in \bR^{n_2 \times r}$.
The tangent space of $\cM_r$ at $\m{X}$ is
\equ{
	T_{\m{X}}
	:= T_{\m{X}} \cM_r = \left\{
	\m{U}\m{B}^T+\m{C}\m{V}^T:
	\m{B}\in\bR^{n_2\times r},
	\;
	\m{C}\in\bR^{n_1\times r}
	\right\}. \label{eq:Tangentspace}
}
The orthogonal projection onto $T_{\m{X}}$ is given by
\equs{
	\cP_{T_{\m{X}}} (\m{Z})
	=
	\m{U}\m{U}^T\m{Z}
	+
	\m{Z}\m{V}\m{V}^T
	-
	\m{U}\m{U}^T\m{Z}\m{V}\m{V}^T.
} 

\subsubsection{Riemannian gradient}
\noindent

Let $f: \bR^{n_1\times n_2} \to \bR$ be a smooth function. For $\m{\eta} \in T_{\m{X}}$, let $\text{D} f(\m{X})[\m{\eta}]$ denote the directional derivative of $f$ at $\m{X}$ along $\m{\eta}$.
According to \cite[equation (3.31)]{absil2008optimization}, the Riemannian gradient $\operatorname{grad}f(\m{X})\in T_{\m{X}}$ is defined by 
\equs{
	\left\langle
	\operatorname{grad}f(\m{X}),
	\m{\eta}
	\right\rangle = \text{D} f(\m{X})[\m{\eta}],
	\qquad
	\forall \m{\eta}\in T_{\m{X}}.
}
According to \cite[equation (3.37)]{absil2008optimization}, the Riemannian gradient is the tangent-space projection of the Euclidean gradient:
\equs{
	\operatorname{grad}f(\m{X}) = \cP_{T_{\m{X}}} \nabla f(\m{X}).
}

For the matrix sensing objective $f(\m{X}) = \frac{1}{2} \twon{ \cA (\m{X}) - \m{y} }^2$, we have $\nabla f(\m{X}) = \cA^*\sbra{\cA(\m{X})-\m{y}}$ and thus
\equ{
	\operatorname{grad}f(\m{X})
	= \cP_{T_{\m{X}}} \cA^*\sbra{\cA(\m{X})-\m{y}}. \label{eq:gradient}
}

\subsubsection{Riemannian Hessian}
\noindent

According to \cite[Section 2]{absil2013extrinsic}, the Riemannian Hessian satisfies
\equs{
	\begin{split}
		\operatorname{Hess}f(\m{X})[\m{\eta}]
		& = \cP_{T_{\m{X}}}
		\text{D} \sbra{\operatorname{grad}f}(\m{X})[\m{\eta}] \\
		& = \cP_{T_{\m{X}}}
		\text{D} (\cP_{T_{\m{X}}} \nabla f) (\m{X}) [\m{\eta}] \\
		& = \cP_{T_{\m{X}}}
		(\text{D} \cP_{T_{\m{X}}}[\m{\eta}] )\,\nabla f(\m{X}) 
		+ \cP_{T_{\m{X}}} \nabla^2 f(\m{X})[\m{\eta}].
	\end{split}
}
Here, $\text{D} \cP_{T_{\m{X}}}[\m{\eta}]$ is  the directional derivative of 
$P_{T_{\m{X}}}$ along $\m{\eta}$ and $\nabla^2 f(\m{X})[\m{\eta}] := \text{D} ( \nabla f) (\m{X}) [\m{\eta}]$ denotes the Euclidean Hessian applied to the direction $\m{\eta}$.

For the matrix sensing objective, we have $\nabla^2 f(\m{X})[\m{\eta}] = \cA^*\cA(\m{\eta})$ and thus
\equ{
		\operatorname{Hess}f(\m{X})[\m{\eta}]
		= \cP_{T_{\m{X}}}
		(\text{D} \cP_{T_{\m{X}}}[\m{\eta}]) \, 
		\cA^*\sbra{\cA(\m{X})-\m{y}} 
		+ \cP_{T_{\m{X}}} \cA^*\cA(\m{\eta}). \label{eq:Hessian}
}

\subsection{RGN algorithm}
The problem in \eqref{eq:main_problem} can be equivalently formulated as the following optimization problem over the fixed-rank matrix manifold:
\equ{
	\min_{\m{X} \in \cM_r }
	\frac{1}{2}
	\twon{ \cA (\m{X}) - \m{y} }^2.
	\label{eq:main_problem1}
}

Motivated by the local quadratic convergence of Newton-type methods, we first consider the Riemannian Newton method for \eqref{eq:main_problem1}. Specifically, at the $t$-th iterate $\m{X}_t\in\cM_r$, let $T_t:=T_{\m{X}_t}$. According to \cite[Section 6.2]{absil2008optimization}, the Riemannian Newton direction $\m{\eta}_t\in T_t$ satisfies 
\equ{
	\operatorname{Hess}f(\m{X}_t)[\m{\eta}_t]
	= -\operatorname{grad}f(\m{X}_t). \label{eq:Newton}
}
Substituting \eqref{eq:gradient} and \eqref{eq:Hessian} into \eqref{eq:Newton} yields
\equs{
	\cP_{T_t} (\text{D} \cP_{T_t}[\m{\eta}_t]) \, \cA^*\sbra{\cA(\m{X}_t)-\m{y}} + \cP_{T_t}\cA^*\cA(\m{\eta}_t) = -\cP_{T_t}\cA^*
	\sbra{\cA(\m{X}_t)-\m{y}}.
}
The first term is the curvature contribution caused by the variation of the tangent space. It is linear in the residual $\cA(\m{X}_t)-\m{y}$ and vanishes at $\m{X}_\star$. To reduce the computational complexity, the RGN approximation discards this curvature term and considers
\equ{
	\cP_{T_t}\cA^*\cA (\m{\eta}_t )
	= -\cP_{T_t}\cA^*
	\sbra{\cA(\m{X}_t)-\m{y}}, \qquad \m{\eta}_t \in T_t. \label{eq:normal_equation0}
}

Define the variable $\m{Z} := \m{X}_t + \m{\eta}_t$. It follows that $\m{Z} \in T_t$. Then, \eqref{eq:normal_equation0} can be rewritten as the RGN equation
\equ{
	\cP_{T_t}\cA^*\cA \cP_{T_t} ( \m{Z} )
	= \cP_{T_t}\cA^* \m{y}, \qquad \m{Z} \in T_t. \label{eq:normal_equation}
}
We let $\m{W}_{t+1}$ denote the solution to this linear equation.

Since every matrix in $T_t$ has rank at most $2r$, the update
$\m{W}_{t+1}$ generally does not belong to $\cM_r$.
To map $\m{W}_{t+1}$ back onto $\cM_r$, we employ the best rank-$r$
approximation as a retraction. Specifically, let $\cT_r$ denote the
best rank-$r$ approximation obtained via the SVD. We set
\equs{
	\m{X}_{t+1}=\cT_r(\m{W}_{t+1}).
}

The complete RGN algorithm with spectral initialization for low-rank matrix recovery is summarized
in Algorithm \ref{alg:RGN}.

\begin{algorithm}
	\caption{RGN algorithm for low-rank matrix recovery}
	\label{alg:RGN}
	\begin{algorithmic}[1]
		\REQUIRE Observations $\m{y}$, sensing operator $\cA$, rank $r$ \\
		\STATE Initialize $\m{X}_0 = \cT_r \sbra{ \cA^* \m{y} }$
		\FOR{$t = 0,1,...$}
		\STATE Obtain the solution $\m{W}_{t+1}$ to \eqref{eq:normal_equation}
		\STATE Set $\m{X}_{t+1} = \cT_r(\m{W}_{t+1})$
		\ENDFOR
		\ENSURE $\m{X}_{t+1}$
	\end{algorithmic}
\end{algorithm}

\begin{remark}[Connections with RGD]
	Compared with RGN, RGD \cite{wei2016guarantees} differs only in the computation of the intermediate iterate $\m{W}_{t+1}$. Specifically, RGD updates
	\equs{
		\begin{split}
			\m{W}_{t+1}
			&= \m{X}_t-\alpha_t\operatorname{grad}f(\m{X}_t)\\
			&= \m{X}_t-\alpha_t\cP_{T_t}\cA^*
			\sbra{\cA(\m{X}_t)-\m{y}},
		\end{split}
	}
	where $\alpha_t$ is the stepsize. Since the search direction is determined solely by the Riemannian gradient, RGD is a first-order method and exhibits linear convergence.
	
	In contrast, RGN is derived from the Riemannian Newton equation and incorporates second-order information through the Gauss-Newton approximation to the Riemannian Hessian. This enables RGN to achieve quadratic local convergence, as established in Section \ref{sec:theory}.
\end{remark}

\section{Theoretical Guarantee} \label{sec:theory}

\subsection{Main Result}
We have the following recovery guarantee for RGN.
\begin{theorem} \label{thm:main}
	Let $\m{X}_\star \in \bR^{n_1 \times n_2}$ be a rank-$r$ matrix with condition number $\kappa$, and let $\m{y}=\cA (\m{X}_\star) \in \bR^m$, where $\cA$ is defined in \eqref{eq:sensing} and the entries of $\{ \m{A}_i \}$ are i.i.d. $\cN(0,1/m)$. Let $\{\m{X}_t\}_{t \ge 0}$ be the RGN iterates generated by Algorithm \ref{alg:RGN}.
	There exists a universal constant $C > 0$ such that, if
	\equs{
		m \ge C r(n_1+n_2) \kappa^2 ,
	}
	then, with probability at least $1-11 \exp\sbra{-(n_1+n_2)/2}$, the iterates satisfy
	\equs{
		 \frobn{\m{X}_t - \m{X}_\star} \le \sqrt{2r} \, 2^{-2^t} \sigma_r(\m{X}_\star), \qquad \forall t \ge 0.
	}
\end{theorem}
\begin{proof}
	The proof is deferred to Section \ref{subsec:proof}.
\end{proof}

Theorem \ref{thm:main} establishes that RGN simultaneously achieves rank-optimal sample complexity and quadratic convergence for Gaussian matrix sensing.

\emph{Rank-optimal sample complexity.}
A rank-$r$ matrix in $\bR^{n_1\times n_2}$ has $r(n_1+n_2-r)$ degrees of freedom, which gives an information-theoretic lower bound on the number of samples required for recovery. The sample size required by RGN matches this lower bound up to the factor $\kappa^2$. 

\emph{Quadratic convergence.}
Achieving rank-optimal sample complexity does not come at the expense of the convergence rate. Specifically, the recovery error of RGN decays doubly exponentially as $2^{-2^t}$, so only $\cO\sbra{\log\log(1/\varepsilon)}$ iterations are required to achieve
$\frobn{\m{X}_t-\m{X}_\star} \le \varepsilon\sqrt{r} \sigma_r(\m{X}_\star)$,
which implies
$\frobn{\m{X}_t-\m{X}_\star} \le \varepsilon\frobn{\m{X}_\star}$.

\subsection{Preliminaries for the Proof}
Restricted isometry property (RIP) plays a key role in controlling the sensing operator in low-rank matrix recovery.
\begin{definition}[RIP] \label{def:RIP}
	The rank-$r$ isometry constant of a linear operator
	$\cA:\bR^{n_1\times n_2}\to\bR^m$
	is defined as the smallest constant $\delta_r\ge0$ such that
	\equs{
		(1-\delta_r)\frobn{\m Z}^2
		\le
		\twon{\cA(\m Z)}^2
		\le
		(1+\delta_r)\frobn{\m Z}^2
	}
	for every $\m Z\in\bR^{n_1\times n_2}$ with $\rank(\m Z)\le r$.
	The operator $\cA$ is said to satisfy the rank-$r$ RIP if $\delta_r<1$.
\end{definition}

The Gaussian sensing operator satisfies the RIP with an optimal sample complexity, as stated below.

\begin{lemma}[\cite{candes2011tight}] \label{lem:gaussian-RIP}
	Define the operator $\cA$ by \eqref{eq:sensing}, where the entries of $\{ \m{A}_i \}$ are i.i.d. $\cN(0,1/m)$.
	For any positive integer $r$ and any $\delta \in (0,1)$, $\cA$ satisfies the rank-$r$ RIP with $\delta_r \le \delta$ with probability at least $1-\exp\sbra{-(n_1+n_2)}$, provided
	\equs{
		m \gtrsim \delta^{-2}r(n_1+n_2).
	}
\end{lemma}

The following lemma establishes restricted-isometry bounds for the tangent-space operator $\cP_T\cA^*\cA\cP_T$ and its inverse.

\begin{lemma} \label{lem:tangent_RIP}
	Let $\cP_T$ be the orthogonal projection onto the tangent space $T$ of the rank-$r$ matrix manifold. Assume that $\cA$ satisfies the rank-$2r$ RIP with constant $\delta_{2r}<1$.
	Then, $\cP_T \cA^*\cA \cP_T : T \to T$ is positive definite, and, for every $\m{Z}\in T$,
	\equs{
		\begin{split}
			(1-\delta_{2r}) \frobn{\m{Z}}
			& \le \frobn{ \cP_T \cA^*\cA \cP_T (\m{Z}) }
			\le (1+\delta_{2r}) \frobn{\m{Z}}, \\
			\frac{1}{1+\delta_{2r}} \frobn{\m{Z}}
			& \le \frobn{ \mbra{\cP_T \cA^*\cA \cP_T}^{-1} (\m{Z}) }
			\le \frac{1}{1-\delta_{2r}} \frobn{\m{Z}}.
		\end{split}
	}
\end{lemma}

\begin{proof}
	The proof is deferred to Appendix \ref{append:tangent_RIP}.
\end{proof}

Lemma \ref{lem:tangent_RIP} shows that, when $\cA$ satisfies the $2r$-RIP, the operator $\cP_{T_t}\cA^*\cA\cP_{T_t}$ is invertible. Consequently, \eqref{eq:normal_equation} admits a unique solution
\equ{
	\m{W}_{t+1} 
	= \mbra{\cP_{T_t} \cA^*\cA \cP_{T_t}}^{-1} \cP_{T_t} \cA^* \m{y}. \label{eq:W0}
}

\subsection{Main Idea of the Proof}
This section outlines the main ideas and motivation underlying the proof of Theorem \ref{thm:main}. 

Under a mild RIP condition, RGN converges quadratically to $\m{X}_\star$ once an iterate enters a sufficiently small Frobenius-norm neighborhood of $\m{X}_\star$, as stated in the following lemma.

\begin{lemma} \label{lem:error_contraction_F}
	Assume that $\cA$ satisfies the rank-$3r$ RIP with constant $\delta_{3r} \le 1/2$. Then, for $t \ge 0$, 
	\equs{
		\frobn{ \m{X}_{t+1} - \m{X}_\star } 
		\le \frac{4}{\sigma_r(\m{X}_\star)} \frobn{ \m{X}_t - \m{X}_\star }^2.
	}
	Further, if
	\equ{
		\frobn{\m{X}_T-\m{X}_\star}
		\le \frac{\sigma_r(\m{X}_\star)}{8} \label{eq:XT_Xstar}
	}
	for some $T\ge0$, then
	\equs{
		\begin{split}
			\frobn{\m{X}_t-\m{X}_\star}
			& \le \frac{\sigma_r(\m{X}_\star)}{4}
			\sbra{
				\frac{4}{\sigma_r(\m{X}_\star)}
				\frobn{\m{X}_T-\m{X}_\star}
			}^{2^{t-T}} \\
			& \le 2^{-2^{t-T}}\sigma_r(\m{X}_\star),
			\qquad \forall t\ge T.
		\end{split}
	}
\end{lemma}
\begin{proof}
	The proof is deferred to Appendix \ref{append:error_contraction_F}.
\end{proof}

By Lemma \ref{lem:gaussian-RIP}, the condition $\delta_{3r} \le 1/2$ in Lemma \ref{lem:error_contraction_F} holds with high probability under the optimal sample complexity $m \gtrsim r(n_1+n_2)$. Therefore, it remains to establish the neighborhood condition \eqref{eq:XT_Xstar} for some $T\ge 0$. 

A direct analysis of the spectral initializer does not yield the desired sample complexity. Indeed, if we set $T=0$ and use the spectral initialization $\m{X}_0 = \cT_r(\cA^*\m{y})$, then it follows from \cite[equation (4.15)]{wei2016guarantees} that $$\frobn{ \m{X}_0 - \m{X}_\star } \le 2 \delta_{2r} \frobn{\m{X}_\star} \le 2 \delta_{2r} \sqrt{r} \kappa \sigma_r(\m{X}_\star).$$ Thus, to guarantee \eqref{eq:XT_Xstar}, one would need $\delta_{2r} \lesssim 1/(\sqrt{r}\kappa)$. By Lemma \ref{lem:gaussian-RIP}, this condition requires $m \gtrsim r^2 (n_1+n_2) \kappa^2$, which is suboptimal because of its quadratic dependence on $r$.

To remove the additional factor of $r$, we instead study the spectral-norm error after $T$ iterations. Suppose that we can establish
\equ{
	\twon{ \m{X}_T - \m{X}_\star } \le 2^{-2^{T}} \sigma_r(\m{X}_\star), \label{eq:XT_Xstar_2}
} 
for $T = \lceil \log_2 \log_2 (16r) \rceil$, 
then \eqref{eq:XT_Xstar} is achieved since
\equ{
	\begin{split}
	\frobn{ \m{X}_T - \m{X}_\star } 
	\le \sqrt{2r} \, \twon{ \m{X}_T - \m{X}_\star }
	\le \sqrt{2r} \, 2^{-2^{T}} \sigma_r(\m{X}_\star)
	& \le \sqrt{2r} (16r)^{-1} \sigma_r(\m{X}_\star) \\
	& \le \sigma_r(\m{X}_\star) / 8. 
	\end{split} \label{eq:XT_Xstar_F}
}

Now, our goal is to show \eqref{eq:XT_Xstar_2}, which is the main technical part of the proof. This is nontrivial for two reasons. First, the desired contraction must be established in the spectral norm rather than the Frobenius norm. Second, the analysis must work under the rank-optimal sample complexity $m \gtrsim r(n_1+n_2) \kappa^2$. 

To show \eqref{eq:XT_Xstar_2}, we first consider establishing a one-step contraction. Given $t>0$, we have
\equ{
	\begin{split}
		& \twon{ \m{X}_{t+1} - \m{X}_\star } \\
		& \le
		\twon{\m{X}_{t+1} - \m{W}_{t+1}}
		+ \twon{\m{W}_{t+1} - \m{X}_{\star}}
		\\
		& \le 2 \twon{\m{W}_{t+1} - \m{X}_{\star}} \\
		& \le 2 \twon{ (\cP_{T_t} - \cI) \m{X}_{\star} } + 2 \twon{ \mbra{\cP_{T_t} \cA^*\cA \cP_{T_t}}^{-1} \cP_{T_t} \sbra{\cA^*\cA - \cI} (\cI - \cP_{T_t}) \m{X}_{\star} },
	\end{split} \label{eq:Xt_Xstar_2}
}
where the second inequality holds since $\rank \sbra{\m{X}_{\star}} = r$ and $\m{X}_{t+1}$ is a best rank-$r$ approximation of $\m{W}_{t+1}$, and the last holds since it follows from \eqref{eq:W0} and $\m{y} = \cA (\m{X}_{\star})$ that
\equ{
	\m{W}_{t+1} 
	= \cP_{T_t} \m{X}_{\star} + \mbra{\cP_{T_t} \cA^*\cA \cP_{T_t}}^{-1}
	\cP_{T_t} \sbra{\cA^*\cA - \cI}
	(\cI - \cP_{T_t}) \m{X}_{\star}.
	\label{eq:W}
}

The first term in the last line of \eqref{eq:Xt_Xstar_2} can be well bounded using Lemma \ref{lem:wei}. The main difficulty lies in controlling the second term, which gives rise to the following two challenges:
\begin{enumerate}[label=(\roman*)]
	\item obtaining a sharp spectral-norm bound for $\twon{ \sbra{\cA^*\cA - \cI} (\cI - \cP_{T_t}) \m{X}_{\star} } $;
	\item controlling the inverse operator $\mbra{\cP_{T_t} \cA^*\cA \cP_{T_t}}^{-1}$ in the spectral norm.
\end{enumerate}

\subsubsection{Virtual-sequence decoupling for Challenge (i)} \label{subsec:virtual}
\noindent

The key to address Challenge (i) is to establish the sharp spectral-norm estimate
\equ{
	\twon{ \sbra{\cA^*\cA-\cI} \m{M} } \lesssim \twon{ \m{M} }, \label{eq:AAI}
}
where $\m{M} := (\cI - \cP_{T_t}) \m{X}_{\star}$ has rank at most $r$. 

A direct uniform concentration argument, such as the RIP, is insufficient for proving \eqref{eq:AAI}. Specifically, Lemma \ref{lem:proj-RIP} yields
\equs{
	\twon{ \sbra{\cA^*\cA-\cI} \m{M} } 
	\le \delta_{r+1} \frobn{ \m{M} } 
	\le \sqrt{r}\,\delta_{r+1} \twon{ \m{M} }.
}
Thus, to guarantee \eqref{eq:AAI}, we require $\delta_{r+1} \lesssim 1/\sqrt{r}$. By Lemma \ref{lem:gaussian-RIP}, this condition translates into $m \gtrsim r^2 (n_1+n_2)$, which is suboptimal in its dependence on $r$. In fact, \cite[Theorem 3.1]{stoger2024non} provides a negative result showing that, under the optimal sample complexity $m \gtrsim r (n_1+n_2)$, \eqref{eq:AAI} cannot,  in general, hold uniformly over all matrices of rank at most $r$.

Note that the matrix $\m{M} = (\cI - \cP_{T_t}) \m{X}_{\star}$ in our problem is not an arbitrary low-rank matrix. Instead, it is determined by the RGN iterate $\m{X}_t$ and therefore belongs to a trajectory-dependent class. Consequently, it suffices to establish \eqref{eq:AAI} along the RGN trajectory, rather than uniformly over the entire set of low-rank matrices.

To exploit this additional structure, we estimate the spectral norm via an $\varepsilon$-net argument. Specifically, Let $\cN_1$ and $\cN_2$ be $1/4$-nets of $\bS^{n_1-1}$ and $\bS^{n_2-1}$, respectively, and set $\cN := \cN_1 \times \cN_2$. By the standard $\varepsilon$-net bound for the spectral norm \cite[Exercise 4.4.3]{vershynin2018high}, we have
\equ{
	\begin{split}
		\twon{ (\cA^*\cA-\cI) \m{M} } 
		& \le 2 \sup_{(\m{x},\m{y}) \in \cN} \abs{ \langle \m{x}\m{y}^T, (\cA^*\cA-\cI) \m{M} \rangle } \\
		& = 2 \sup_{(\m{x},\m{y}) \in \cN} \left|  \sum^m_{i=1} \langle \m{A}_i, \m{M} \rangle \langle \m{x}\m{y}^T,  \m{A}_i \rangle - \langle \m{x}\m{y}^T, \m{M} \rangle \right|.
	\end{split} \label{eq:net}
}
Note that $\m{M}$ is determined by $\m{X}_t$ and hence depends on $\{ \langle \m{x}\m{y}^T,  \m{A}_i \rangle \}$. This dependence prevents an application of Gaussian concentration inequalities to \eqref{eq:net}.

To remove this dependence, inspired by \cite{stoger2024non,cai2025fast}, we exploit the rotation invariance of Gaussian random variables and construct an auxiliary virtual sequence. 
Specifically, for fixed $(\m{x},\m{y})\in\cN$ and any $\m{Z}\in \bR^{n_1\times n_2}$, we define the orthogonal projection operators
\equs{
	\tilde{\cP}\m{Z}:=\langle \m{x}\m{y}^T,\m{Z}\rangle\m{x}\m{y}^T, \qquad \qquad \tilde{\cP}^{\perp}:=\cI-\tilde{\cP}.
}
Here and throughout this paper, the tilde indicates dependence on the fixed pair $(\m{x},\m{y})$, which is
suppressed for notational simplicity.

Since the entries of $\m{A}_i$ are i.i.d. Gaussian, it follows from the rotation invariance property that $\{ \tilde{\cP}^\perp \m{A}_i \}$ are stochastically independent of $\{ \langle \m{x}\m{y}^T,  \m{A}_i \rangle \}$ \cite{stoger2024non,cai2025fast}. 
This fact motivates us to construct a virtual sensing operator $\tilde{\cA}$ using only $\{ \tilde{\cP}^\perp \m{A}_i \}$. Specifically, define $\tilde{\cA} : \bR^{n_1 \times n_2} \to \bR^{m+1}$ by
\equ{
	[ \tilde{\cA} (\m{Z}) ]_i =
	\begin{cases}
		\displaystyle \langle \tilde{\cP}^{\perp} \m{A}_i,\m{Z} \rangle, & i=1,\ldots,m,\\[1ex]
		\langle \m{x}\m{y}^T,\m{Z} \rangle, & i=m+1.
	\end{cases}
	\label{eq:virtual-sensing-operator}
}
The additional $(m+1)$-st coordinate compensates for the rank-one direction removed by $\tilde{\cP}^{\perp}$ in expectation, ensuring that $\tilde{\cA}$ remains isotropic. Indeed, as $\bE[\cA^* \cA] = \cI$, we have $
\bE [ \tilde{\cA}^* \tilde{\cA} ] = \bE [ \tilde{\cP} + \tilde{\cP}^{\perp} \cA^* \cA \tilde{\cP}^{\perp} ] = \cI$. 

We then generate a virtual RGN sequence $\{ \tilde{\m{X}}_t \}$ by replacing $\cA$ with $\tilde{\cA}$ in Algorithm \ref{alg:RGN}; see Algorithm \ref{alg:virtual}. For each fixed $(\m{x},\m{y})$, $\{ \tilde{\m{X}}_t \}$ depends only on $\{ \tilde{\cP}^\perp \m{A}_i \}$ and $\m{X}_\star$, and is therefore independent of $\{ \langle \m{x}\m{y}^T,  \m{A}_i \rangle \}$. The same independence holds for the corresponding virtual matrix $\tilde{\m{M}} := (\cI - \cP_{\tilde{T}_t}) \m{X}_{\star}$.

\begin{algorithm}[htb]
	\caption{Virtual RGN sequence for a fixed $(\m{x},\m{y}) \in \cN$} \label{alg:virtual}
	\begin{algorithmic}[1]
		\STATE Initialize $\tilde{\m{X}}_0 = \cT_r ( \tilde{\cA}^*\tilde{\cA} (\m{X}_\star) )$
		\FOR{$t = 0,1,...$}
		\STATE Let $\tilde{T}_t$ be the tangent space of $\tilde{\m{X}}_t$. Obtain the solution $\tilde{\m{W}}_{t+1}$ to the equation 
		\equs{
			\cP_{\tilde{T}_t} \tilde{\cA}^* \tilde{\cA} \cP_{\tilde{T}_t} ( \m{Z} )
			= \cP_{\tilde{T}_t} \tilde{\cA}^* \tilde{\cA} (\m{X}_\star), \qquad \m{Z} \in \tilde{T}_t.
		}
		\vspace{-1em}
		\STATE Set $\tilde{\m{X}}_{t+1} = \cT_r ( \tilde{\m{W}}_{t+1} )$
		\ENDFOR
	\end{algorithmic}
\end{algorithm}

Using this virtual-sequence, we can sharply bound \eqref{eq:net}. Specifically, Lemma \ref{lem:spectral-control} shows that with probability at least $1- 3 \exp(-(n_1+n_2))$, we have
\equ{
	\twon{(\cA^*\cA-\cI)\m{M}}
	\lesssim
	\twon{\m{M}}
	+
	\sup_{(\m{x},\m{y})\in\cN}
	\frobn{ \m{M}-\tilde{\m{M}} }, \label{eq:AAI_simplify}
}
provided $m \gtrsim r(n_1+n_2)$. 
Thus, the desired estimate \eqref{eq:AAI} is recovered up to an additional error term, measuring the distance between $\m{M}$ and $\tilde{\m{M}}$.
This distance is directly related to the perturbation of the tangent spaces $T_t$ at $\m{X}_t$ and $\tilde{T}_t$ at $\tilde{\m{X}}_t$.
This motivates us to consider the joint error 
\equ{
	E_t := \twon{ \m{X}_t - \m{X}_\star } + \sup_{(\m{x},\m{y}) \in \cN} \frobn{ \m{X}_t - \tilde{\m{X}}_t }, \label{eq:Et}
}
rather than the single error $\twon{ \m{X}_t - \m{X}_\star }$, in the subsequent analysis.

\subsubsection{Finite Neumann expansion for Challenge (ii)}
\noindent

For Challenge (ii), we point out that the Frobenius-to-Frobenius bound in Lemma \ref{lem:tangent_RIP} for the inverse operator is not sufficiently sharp and will result in a suboptimal rank-quadratic sample complexity. To avoid this loss, we use an exact finite Neumann expansion for the inverse operator, as stated below. 

\begin{lemma}[Finite Neumann expansion] \label{lem:Neumann}
	Let $\cB$ be a linear operator. If $\cI-\cB$ is invertible, then, for every integer $K > 0$,
	\equs{
		(\cI - \cB)^{-1} = \sum^{K-1}_{k=0} \cB^k + \cB^K (\cI - \cB)^{-1}.
	}
\end{lemma} 
\begin{proof}
	We have 
	$(\cI + \cB + \ldots + \cB^{K-1}) ( \cI - \cB ) = \cI - \cB^K$.
	Multiplying both sides by $(\mathcal{I}-\mathcal{B})^{-1}$
	yields the desired identity.
\end{proof}

As shown in Lemma \ref{lem:tangent_RIP}, $\cP_{T_t} \cA^*\cA \cP_{T_t}$ is an invertible operator on $T_t$. Thus, we have
\equs{
	\mbra{\cP_{T_t} \cA^*\cA \cP_{T_t}}^{-1} = \mbra{ \cI - \cP_{T_t} \sbra{\cI - \cA^*\cA} \cP_{T_t} }^{-1}.
}
It then follows from Lemma \ref{lem:Neumann} that
\equ{
	\begin{split}
		\mbra{\cP_{T_t} \cA^*\cA \cP_{T_t}}^{-1}
		= \sum_{k=0}^{K-1}
		\mbra{\cP_{T_t} \sbra{\cI - \cA^*\cA} \cP_{T_t}}^k + \mbra{\cP_{T_t} \sbra{\cI - \cA^*\cA} \cP_{T_t}}^K
		\mbra{\cP_{T_t} \cA^*\cA \cP_{T_t}}^{-1}.
	\end{split}
	\label{eq:neumann}
}

For a fixed iteration $t$, we define
\equ{
	\begin{split}
		\m{Z}_0 & := \cP_{T_t} \sbra{\cA^*\cA - \cI} (\cI - \cP_{T_t}) \m{X}_{\star}, \\
		\m{Z}_{k+1} & := \cP_{T_t} \sbra{\cI - \cA^*\cA} \m{Z}_k, \qquad k \ge 0.
	\end{split} \label{eq:Z0Zk}
}
Using \eqref{eq:neumann} and \eqref{eq:Z0Zk} yield
\equ{
	\begin{split}
		& \twon{ \mbra{\cP_{T_t} \cA^*\cA \cP_{T_t}}^{-1} \cP_{T_t} \sbra{\cA^*\cA - \cI} (\cI - \cP_{T_t}) \m{X}_\star } \\
		& \le \sum_{k=0}^{K-1} \twon{\m{Z}_k}  + \twon{ [\cP_{T_t} \sbra{\cI - \cA^*\cA} \cP_{T_t} ]^K \mbra{\cP_{T_t} \cA^*\cA \cP_{T_t}}^{-1} \m{Z}_0 }.
	\end{split} \label{eq:bound1}
}
Here, the finite sum is inverse-free, while the remainder still contains the inverse but is multiplied by a $K$-fold contraction factor $[ \cP_{T_t} \sbra{\cI - \cA^*\cA} \cP_{T_t} ]^K$, ensuring that it can be bounded sharply. We next show how to bound these two terms separately. 

First, we consider the finite sum in \eqref{eq:bound1}. To bound $\twon{\m{Z}_k}$, we need to apply the virtual-sequence argument in addressing Challenge (i). For every $(\m{x},\m{y}) \in \cN$, we define the corresponding virtual recursion by
\equ{
	\begin{split}
		\tilde{\m{Z}}_0 & := \cP_{\tilde{T}_t} (\tilde{\cA}^*\tilde{\cA} - \cI) (\cI - \cP_{\tilde{T}_t}) \m{X}_{\star}, \\
		\tilde{\m{Z}}_{k+1} & := \cP_{\tilde{T}_t} (\cI - \tilde{\cA}^*\tilde{\cA}) \tilde{\m{Z}}_k, \qquad k \ge 0.
	\end{split} \label{eq:Z0Zk_tilde}
}
To bound $\twon{\m{Z}_k}$ by recursively using the trajectory-dependent estimate \eqref{eq:AAI_simplify}, we need to control the difference between $\m{Z}_{k-1}$ and $\tilde{\m{Z}}_{k-1}$. The following lemma provides the required one-step estimate.

\begin{lemma} \label{lem:ZminusZtilde}
	Let $\m{M}$ be a matrix derived from the RGN iterates with $\rank(\m{M}) \le 2r$.
	For each $(\m{x},\m{y})\in\cN$, let
	$\tilde{\m{M}}$ denote the matrix constructed from the associated virtual sequence and satisfying $\rank(\tilde{\m{M}})\le 2r$.
	Define
	\equs{
		\m{Z}
		:=
		\cP_{T_t}(\cA^*\cA-\cI)\m{M},
		\qquad
		\tilde{\m{Z}}
		:=
		\cP_{\tilde{T}_t}
		(\tilde{\cA}^*\tilde{\cA}-\cI)
		\tilde{\m{M}}.
	}
	Suppose that $\delta_{6r} \le 1$ and $E_t \le c_0 \sigma_r(\m{X}_\star)$ for some $0<c_0<1$. Then, with probability at least $1-4\exp(-5(n_1+n_2))$, 
	\equs{
		\sup_{(\m{x},\m{y})\in\cN}
		\frobn{\m{Z}-\tilde{\m{Z}}} 
		\le \frac{(8 c_0+1)\alpha_1}{1-c_0} \Big( \twon{\m{M}} + \sup_{(\m{x},\m{y})\in\cN} \frobn{\m{M}-\tilde{\m{M}}} \Big),
	}
	where $\alpha_1:=16\sqrt{2r(n_1+n_2)/m}+2\delta_{6r}$.
\end{lemma}

\begin{proof}
	The proof is deferred to Appendix \ref{append:ZminusZtilde}.
\end{proof}

Since the finite sum has only $K$ terms, the estimates in \eqref{eq:AAI_simplify} and Lemma \ref{lem:ZminusZtilde} are applied only $K$ times. A union bound over the $K$ recursion levels therefore preserves the desired high-probability guarantee.

Second, we control the remainder in \eqref{eq:bound1}. Due to the presence of the inverse operator, we relax the spectral norm to the Frobenius norm and apply the tangent-space RIP in Lemma \ref{lem:tangent_RIP} and Lemma \ref{lem:PT_I_AA_PT}. Thanks to the $K$-fold contraction factor, the additional $\sqrt{r}$ factor introduced by the relaxation will be absorbed, thereby achieving the rank-optimal sample complexity.

\subsection{Proof of Theorem \ref{thm:main}} \label{subsec:proof}
The proof consists of two stages. In the first stage, we address the two aforementioned challenges and establish quadratic contraction of the joint error $E_t$ in \eqref{eq:Et} during the first $T$ iterations. This yields the desired spectral-norm
estimate in \eqref{eq:XT_Xstar_2} and shows that the iterate $\m{X}_T$ enters the Frobenius-norm local region as in \eqref{eq:XT_Xstar_F}. In the second stage, we invoke Lemma \ref{lem:error_contraction_F} to establish quadratic convergence of RGN
to $\m{X}_\star$ for all $t\ge T$.

We first recall the initialization guarantee.

\begin{lemma}[{\cite[Lemma 6]{cai2025fast}}] \label{lem:error_init}
	Let $0<c_0<10^{-4}$. There exists a universal constant $C>0$ such that, if $m \ge C r(n_1+n_2) \kappa^2$, then, with probability at least $1-4\exp\sbra{-(n_1+n_2)}$,
	\equs{
		E_0 \le c_0 \sigma_r(\m{X}_\star).
	}
\end{lemma}

The following lemma gives the one-step contraction of the joint error.
\begin{lemma} \label{lem:error_contraction}
	Let $0<c_0<10^{-4}$ and fix $t\ge0$. There exists a universal constant $C>0$ such that, if $m \ge C r(n_1+n_2)$, then, with probability at least $1-7 \exp\sbra{-(n_1+n_2)}$, the estimate
	\equs{
		E_{t+1} 
		\le \frac{572}{\sigma_r(\m{X}_\star)} E_t^2
	}
	holds whenever $E_t \le c_0 \sigma_r(\m{X}_\star)$.
\end{lemma}
\begin{proof}
	The proof is deferred to Appendix \ref{append:error_contraction}.
\end{proof}

\textbf{Case 1: $0\le t\le T$.}
Let $T=\lceil\log_2\log_2(16r)\rceil$. Since $572c_0<1$, the one-step estimate in Lemma \ref{lem:error_contraction} preserves the condition $E_t\le c_0\sigma_r(\m{X}_\star)$. Combining Lemma \ref{lem:error_init} and Lemma \ref{lem:error_contraction} and applying induction, we obtain
\equs{
	E_t
	\le \frac{\sigma_r(\m{X}_\star)}{572}
	\sbra{
		\frac{572E_0}{\sigma_r(\m{X}_\star)}
	}^{2^t}
	\le 2^{- 2^{t+1} }\sigma_r(\m{X}_\star),
	\qquad \forall\,0\le t\le T,
}
where the last inequality uses $\frac{572E_0}{\sigma_r(\m{X}_\star)} \le 572\times10^{-4} < 2^{-2}$.
Since Lemma \ref{lem:error_contraction} is invoked for $t=0,\ldots,T-1$, a union bound shows that the preceding estimates fail with probability at most
\equs{
	(4+7T) \exp\sbra{-(n_1+n_2)}
	\le 11 T \exp\sbra{-(n_1+n_2)}
	\le 11 \exp\sbra{-(n_1+n_2)/2},
}
where we use $T\le n_1+n_2\le\exp\sbra{(n_1+n_2)/2}$.

By the definition of $E_t$, we have
\equs{
	\twon{\m{X}_t-\m{X}_\star}
	\le E_t
	\le 2^{- 2^{t+1} }\sigma_r(\m{X}_\star),
	\qquad \forall\,0\le t\le T, 
}
and, since $\rank(\m{X}_t-\m{X}_\star)\le2r$,
\equ{
	\frobn{\m{X}_t-\m{X}_\star}
	\le \sqrt{2r}\twon{\m{X}_t-\m{X}_\star}
	\le \sqrt{2r}\,2^{- 2^{t+1} }\sigma_r(\m{X}_\star),
	\qquad \forall\,0\le t\le T. \label{eq:tF}
}

\textbf{Case 2: $t\ge T$.}
It follows from \eqref{eq:tF} and $T=\lceil \log_2\log_2(16r) \rceil$ that 
\equ{
	\frobn{\m{X}_T - \m{X}_\star} 
	\le \sqrt{2r}\,2^{- 2^{T+1} }\sigma_r(\m{X}_\star) 
	\le \sqrt{2r}\,2^{- 2^{\log_2\log_2(16r)} }\sigma_r(\m{X}_\star) 
	\le \sigma_r(\m{X}_\star) / 8. \label{eq:TF}
}
Then, applying Lemma \ref{lem:error_contraction_F} and \eqref{eq:TF}, we obtain
\equ{
	\begin{split}
		\frobn{\m{X}_t-\m{X}_\star}
		& \le
		\frac{\sigma_r(\m{X}_\star)}{4}
		\sbra{
			\frac{4}{\sigma_r(\m{X}_\star)}
			\frobn{\m{X}_T-\m{X}_\star}
		}^{2^{t-T}} \\
		& \le
		\frac{\sigma_r(\m{X}_\star)}{4}
		\sbra{ 4 \sqrt{2r} \, 2^{-2^{T+1}} }^{2^{t-T}} \\
		& \le
		\frac{\sigma_r(\m{X}_\star)}{4}\,2^{-2^t}
		\le
		\sqrt{2r}\,2^{-2^t}\sigma_r(\m{X}_\star), \qquad \forall\, t \ge T.
	\end{split} \label{eq:tF>T}
}

Combining \eqref{eq:tF} with \eqref{eq:tF>T} yields the desired estimate, completing the proof.

\section{Implementation and Numerical Experiments} \label{sec:experiments}
In this section, we first provide the implementation of RGN and then conduct extensive experiments to evaluate its performance.

\subsection{RGN Implementation}
\subsubsection{Step 1: solving the RGN equation \eqref{eq:normal_equation}}
Using a unique parameterization of the tangent space $T_t$, we reformulate \eqref{eq:normal_equation} as an equivalent matrix-form linear system in the parameter space.
Throughout this subsection, we fix an iteration $t$ and omit the subscript
$t$ for notational simplicity.

\emph{Unique parameterization of the tangent space.}
Following the intrinsic representation of tangent vectors on the fixed-rank manifold in \cite{huang2017intrinsic}, we introduce a unique parameterization of $T$ as follows.
The representation of the tangent space in \eqref{eq:Tangentspace} is not unique. For example, the pairs $(\m{B},\m{C})$ and $(\m{B}-\m{V}, \m{C}+\m{U})$ yield the same tangent vector. 
To eliminate this redundancy, we let $\m{U}_{\perp} \in \bR^{n_1 \times (n_1-r)}$ and $\m{V}_{\perp} \in \bR^{n_2 \times (n_2-r)}$ be orthonormal complements of $\m{U}$ and $\m{V}$, respectively, and define the parameter space $$
	\cE := \lbra{ (\m{D},\m{E}_1,\m{E}_2): \m{D} \in \bR^{r \times r}, \m{E}_i \in \bR^{(n_i-r) \times r}, i = 1,2}.$$
Then, we define the linear operator $\cL: \cE \to T$ by
\equs{
	\cL (\m{D},\m{E}_1,\m{E}_2) := \m{U} \m{D} \m{V}^T + \m{U}_\perp \m{E}_1 \m{V}^T + \m{U} \m{E}_2^T \m{V}_\perp^T.
}
The adjoint operator $\cL^* : \bR^{n_1 \times n_2} \to \cE$ is given by
\equs{
	\cL^*(\m{X})
	=
	( \m{U}^T \m{X} \m{V},\;
	\m{U}_{\perp}^T \m{X} \m{V},\;
	\m{V}_{\perp}^T \m{X}^T \m{U} ).
	\label{eq:L_adjoint}
}
A direct calculation gives
\equs{
	\cL \cL^* = \cP_T, \qquad \qquad \cL^* \cL = \cI.
} 
Thus, $\cL^* \cL = \cI$ implies that $\cL$ is injective and isometric. Moreover, $\text{dim}(\cE) = (n_1+n_2-r)r = \text{dim}(T)$. 
Consequently, $\cL$ is bijective, and every element of $T$ admits a unique representation through $\cL$ with parameters in $\cE$.

\emph{Reformulation of \eqref{eq:normal_equation}.}
Using $\cP_T = \cL \cL^*$ and the injectivity of $\cL$, \eqref{eq:normal_equation} is equivalent to
\equ{
	\cL^* \cA^*\cA \cL \cL^* ( \m{Z} )
	= \cL^* \cA^* \m{y}, \qquad \m{Z} \in T.  \label{eq:normal_equation3}
}
Since $\m{Z}\in T$ and $\cL$ is isometric, there exists a unique triple $(\m{D},\m{E}_1,\m{E}_2)\in\cE$ such that
$\m{Z}=\cL(\m{D},\m{E}_1,\m{E}_2)$. Substituting this representation into \eqref{eq:normal_equation3} yields the parameterized equation
\equs{
	(\cA \cL )^* \cA \cL (\m{D},\m{E}_1,\m{E}_2)  = (\cA \cL )^* \m{y}. \label{eq:normal_equation4}
}
We define the matrix $\m{J}\in\bR^{m \times (n_1+n_2-r)r}$ whose $i$th row is
\equs{
	\m{J}_{i,:}
	=
	\begin{bmatrix}
		\operatorname{vec}(\m{U}^T\m{A}_i\m{V})^T &
		\operatorname{vec}\bigl(\m{U}_\perp^T \m{A}_i\m{V}\bigr)^T &
		\operatorname{vec}\bigl(\m{V}_\perp^T \m{A}_i^T \m{U}\bigr)^T
	\end{bmatrix}, \label{eq:J}
}
and let
$\m{z} := \begin{bmatrix} \operatorname{vec}(\m{D}); \operatorname{vec}(\m{E}_1); \operatorname{vec}(\m{E}_2) \end{bmatrix} \in \bR^{(n_1+n_2-r)r}$.
A direct calculation shows
\equs{
	\m{J}\m{z} = \cA \cL(\m{D},\m{E}_1,\m{E}_2). \label{eq:Jz}
}
Consequently, the RGN equation \eqref{eq:normal_equation} is equivalent to the linear system
\equ{
	\m{J}^T\m{J}\m{z}
	=
	\m{J}^T\m{y},
	\label{eq:J_equation}
}
which can be solved using a standard linear solver.
Let $\hat{\m{z}}$ denote the solution to \eqref{eq:J_equation}, and let $(\hat{\m{D}},\hat{\m{E}}_1,\hat{\m{E}}_2)$ denote the corresponding parameter tuple. The resulting update is then given by
$\hat{\m{W}}_{t+1} = \m{U} \hat{\m{D}} \m{V}^T + \m{U}_\perp \hat{\m{E}}_1 \m{V}^T + \m{U} \hat{\m{E}}_2 \m{V}_\perp^T$.

\subsubsection{Step 2: Efficient truncated SVD}
It remains to compute the truncated SVD
$\m{X}_{t+1}=\cT_r(\hat{\m{W}}_{t+1})$.
Inspired by \cite{wei2016guarantees}, we exploit the tangent-space structure of $\hat{\m{W}}_{t+1}$ to avoid forming it explicitly, reducing the computation to the SVD of a $2r\times 2r$ matrix.

Specifically, we compute compact QR factorizations $\m{U}_\perp \hat{\m{E}}_1 = \m{Q}_1 \m{R}_1$ and $\m{V}_\perp \hat{\m{E}}_2 = \m{Q}_2 \m{R}_2$, where $\m{Q}_i \in \bR^{n_i \times r}$ and $\m{R}_i \in \bR^{ r \times r}$ for $i=1,2$. Then,
\equs{
	\begin{split}
		\hat{\m{W}}_{t+1} 
		& = \m{U} \hat{\m{D}} \m{V}^T + \m{Q}_1 \m{R}_1 \m{V}^T + \m{U} \m{R}_2^T \m{Q}_2^T  \\
		& = \begin{bmatrix} \m{U} & \m{Q}_1 \end{bmatrix}
		\begin{bmatrix} \hat{\m{D}} & \m{R}_2^T \\ \m{R}_1 & \m{0} \end{bmatrix}
		\begin{bmatrix} \m{V} & \m{Q}_2 \end{bmatrix}^T.
	\end{split} \label{eq:Wdecomp}
}
Let $\begin{bmatrix}
	\hat{\m{D}}
	&
	\m{R}_2^T \\
	\m{R}_1
	&
	\m{0}
\end{bmatrix} = \m{U}_m
\m{\Sigma}_m \m{V}_m^T$ be the full SVD of this $2r\times 2r$ matrix. Since
$\begin{bmatrix}\m{U}&\m{Q}_1\end{bmatrix}$ and
$\begin{bmatrix}\m{V}&\m{Q}_2\end{bmatrix}$
have orthonormal columns, an SVD of $\hat{\m{W}}_{t+1}$ is
\equs{
	\hat{\m{W}}_{t+1}
	=
	\left(
	\begin{bmatrix}\m{U}&\m{Q}_1\end{bmatrix}\m{U}_m
	\right)
	\m{\Sigma}_m
	\left(
	\begin{bmatrix}\m{V}&\m{Q}_2\end{bmatrix}\m{V}_m
	\right)^T.
}
Thus, $\m{X}_{t+1}$ is obtained by retaining the leading $r$ singular triplets.

\subsubsection{Computational complexity}
For Step 1, constructing $\m{J}\in\bR^{m\times (n_1+n_2-r)r}$ requires $\cO(mn_1n_2r)$ flops. Solving the linear system \eqref{eq:J_equation} using a standard direct solver requires $\cO\bigl(m(n_1+n_2)^2r^2+(n_1+n_2)^3r^3\bigr)$ flops. 
For Step 2, the two compact QR factorizations require $\cO((n_1+n_2)r^2)$ flops, and the SVD of the $2r\times2r$ middle matrix requires $\cO(r^3)$ flops. Recall $n =\max\lbra{n_1,n_2}$. Then, the per-iteration computational cost of RGN is $\cO(mn^2r^2+n^3r^3)$.

\subsection{Numerical Results}
We compare RGN with RISRO \cite{luo2024recursive}, GNMR \cite{zilber2022gnmr}, RGD \cite{wei2016guarantees}, and ScaledGD \cite{tong2021accelerating}. Our MATLAB code for RGN is available at
\begin{center}
	\url{https://github.com/wxm25/RGN}
\end{center}
 We use the publicly available codes of RISRO, GNMR, and ScaledGD, and implement RGD ourselves.

\subsubsection{Phase transition}
To study the sample complexity, we set \(n_1=n_2=50\) and fix the condition number of the ground-truth matrix $\m{X}_\star$ to \(\kappa=2\). The rank varies over \(r\in\{2,4,\ldots,20\}\), while the sample size varies over \(m\in\{400,600,\ldots,2000\}\). The sensing matrix has i.i.d. Gaussian entries distributed as \(\mathcal{N}(0,1/m)\). A trial is declared successful if the relative error satisfies $\frobn{\hat{\m{X}} - \m{X}_\star} / \frobn{\m{X}_\star} <10^{-6}$, where $\hat{\m{X}}$ is the recovered matrix. For each pair $(m, r)$, the success rate is averaged over $20$ independent Monte Carlo trials. 

\begin{figure}
	\centering
		\begin{minipage}{0.31\textwidth}
		\centering
		\includegraphics[height=3cm]{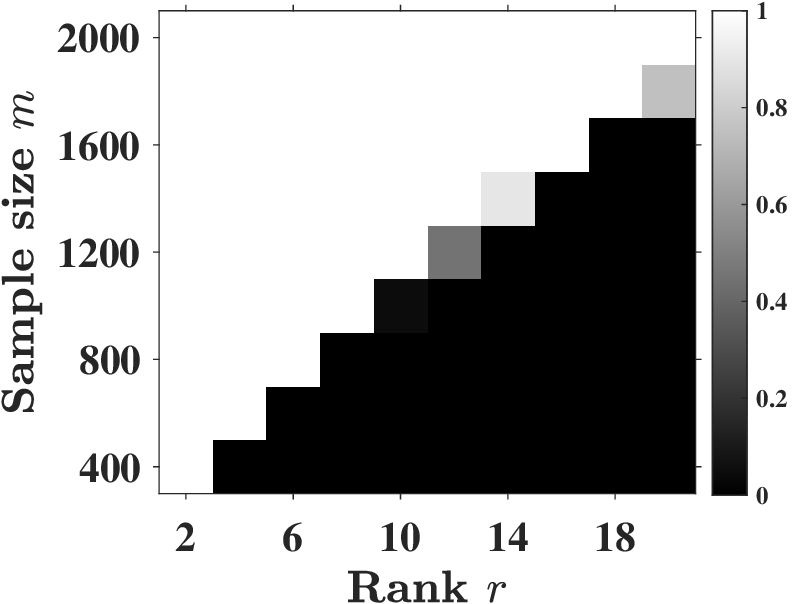}\\
		(a) RGN
	\end{minipage}
	\hfill
	\begin{minipage}{0.31\textwidth}
		\centering
		\includegraphics[height=3cm]{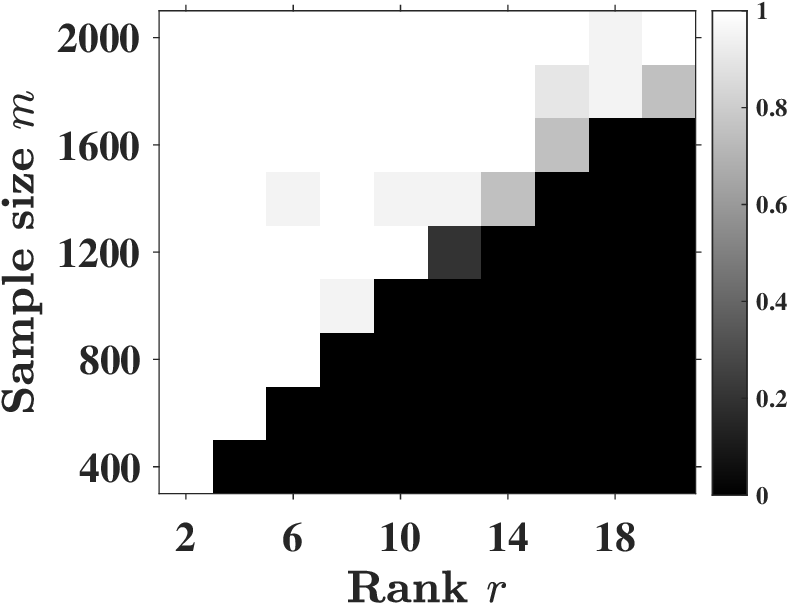}\\
		(b) RISRO \cite{luo2024recursive}
	\end{minipage}
	\hfill
	\begin{minipage}{0.31\textwidth}
		\centering
		\includegraphics[height=3cm]{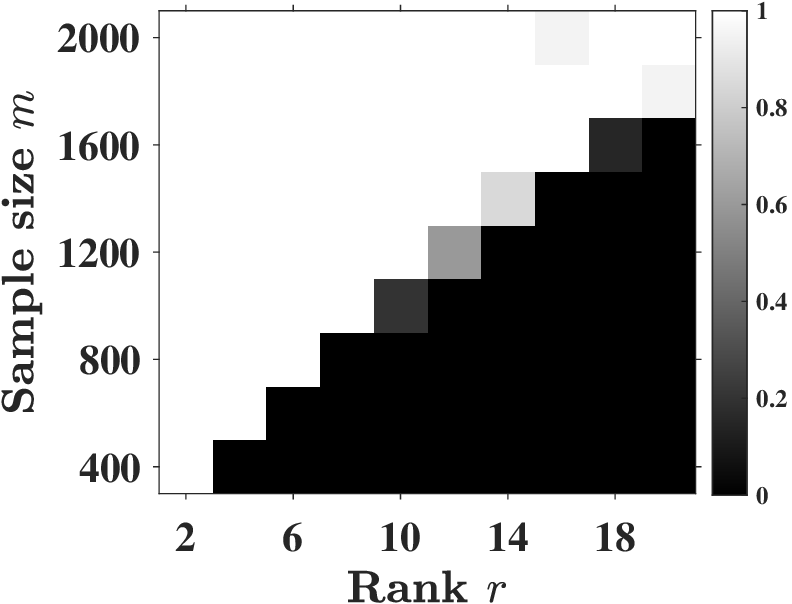}\\
		(c) GNMR \cite{zilber2022gnmr}
	\end{minipage}

	\vspace{0.3cm}
	\begin{minipage}{0.31\textwidth}
		\centering
		\includegraphics[height=3cm]{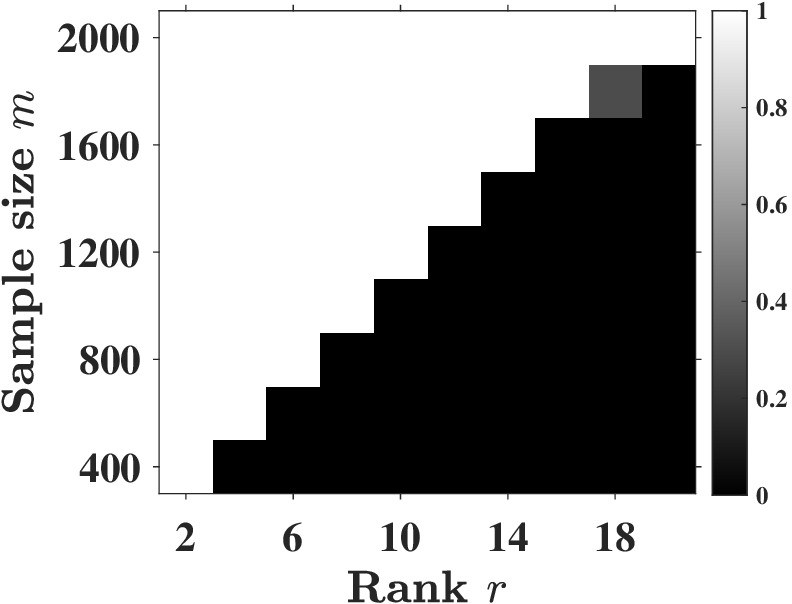}\\
		(d) RGD \cite{wei2016guarantees}
	\end{minipage}
	\hfill
	\begin{minipage}{0.31\textwidth}
		\centering
		\includegraphics[height=3cm]{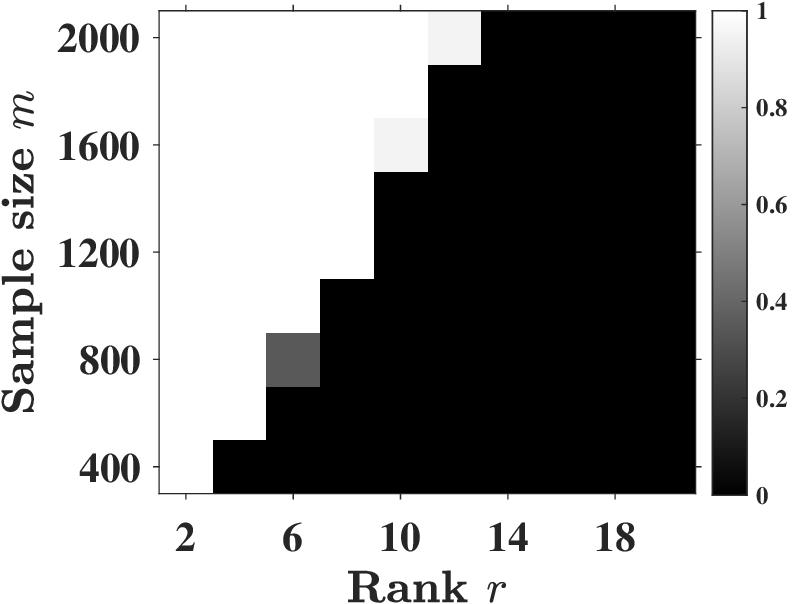}\\
		(e) ScaledGD \cite{tong2021accelerating}
	\end{minipage}
	\hfill
	\begin{minipage}{0.31\textwidth}
		\centering
		\includegraphics[width =3.9cm,height=3cm]{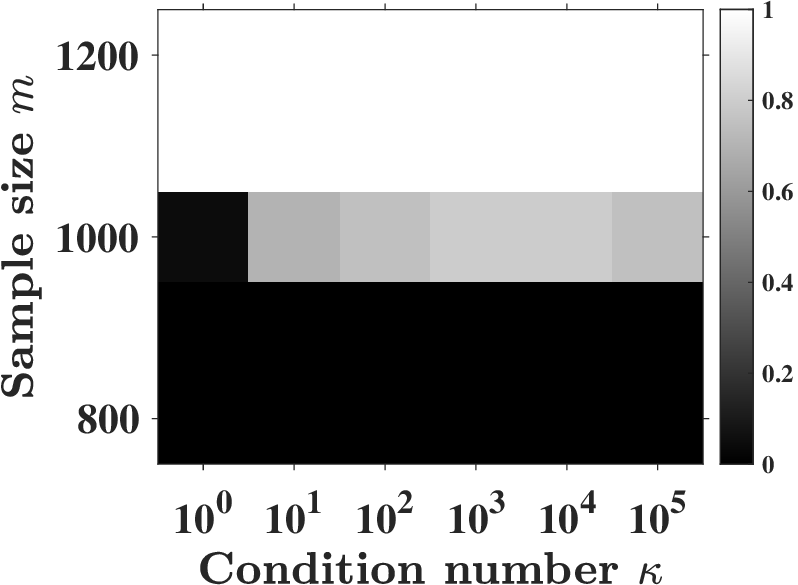}\\
		(f) RGN: varying $\kappa$
	\end{minipage}
	
	\caption{Phase transition results. While means success, while black means failure.}
	\label{fig:phase_transition}
\end{figure}

Figure \ref{fig:phase_transition} illustrates the successful recovery rates. As shown in Figure \ref{fig:phase_transition}(a), the phase transition of RGN exhibits an approximately linear relationship between $m$ and $r$, which is consistent with our theoretical rank-optimal sample complexity. Moreover, RGN achieves better recovery performance than the other methods.

\subsubsection{Convergence performance}
\begin{figure}
	\centering
	\begin{minipage}{0.48\textwidth}
		\centering
		\includegraphics[height=4.5cm]{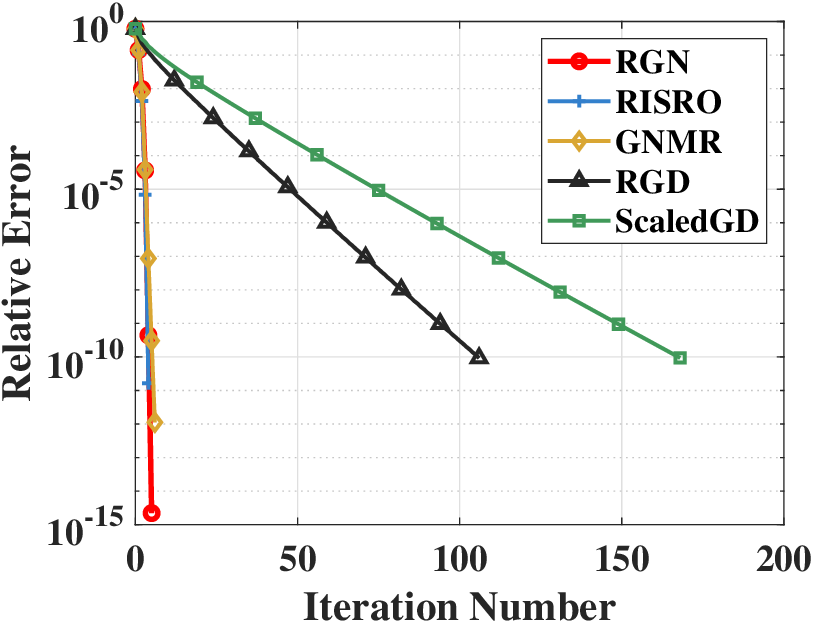}
	\end{minipage}
	\hfill
	\begin{minipage}{0.48\textwidth}
		\centering
		\includegraphics[height=4.5cm]{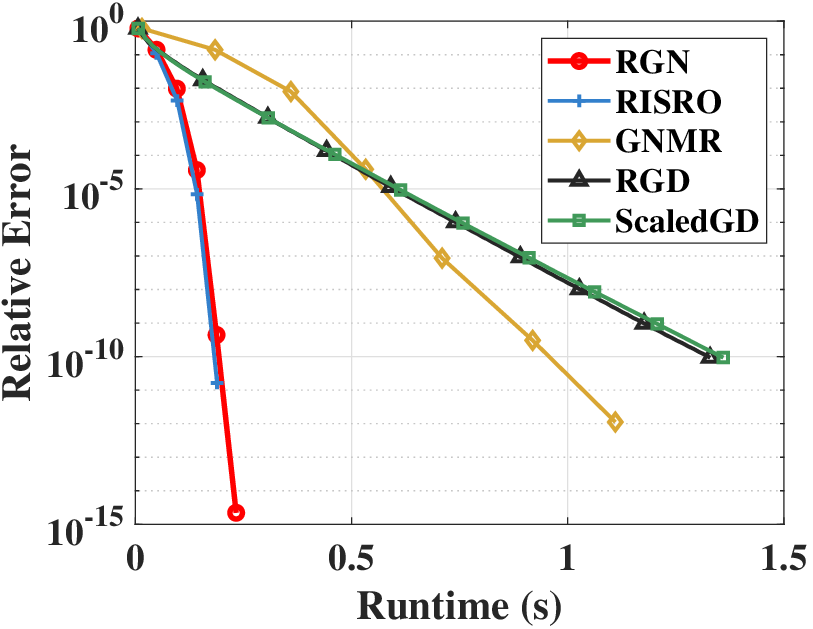}
	\end{minipage}
	\caption{Convergence comparisons.}
	\label{fig:convergence}
\end{figure}

We set $n_1=n_2=100$, $r=5$, $\kappa=2$, and $m=\left\lceil 3.5(n_1+n_2-r)r \right\rceil$.
Figure \ref{fig:convergence} compares the convergence performance of the methods in terms of iteration number and runtime. As shown in Figure \ref{fig:convergence}(a), RGN, RISRO, and GNMR reach high accuracy within only a few iterations, whereas RGD and ScaledGD require substantially more iterations. This observation is consistent with the quadratic convergence of RGN, RISRO, and GNMR, in contrast to the linear convergence of the first-order methods RGD and ScaledGD. Figure \ref{fig:convergence}(b) further shows that RGN and RISRO are faster than other methods.

\subsubsection{Scalability}
We next study the scalability of the methods as the problem size increases. For a fair comparison, the runtime of each method is measured as the time required for successful recovery, and the results are averaged over $20$ trials. We set $\kappa=2$ and $m=\left\lceil 3.5(n_1+n_2-r)r \right\rceil$.

\begin{figure}
	\centering
	\begin{minipage}{0.48\textwidth}
		\centering
		\includegraphics[height=4.5cm]{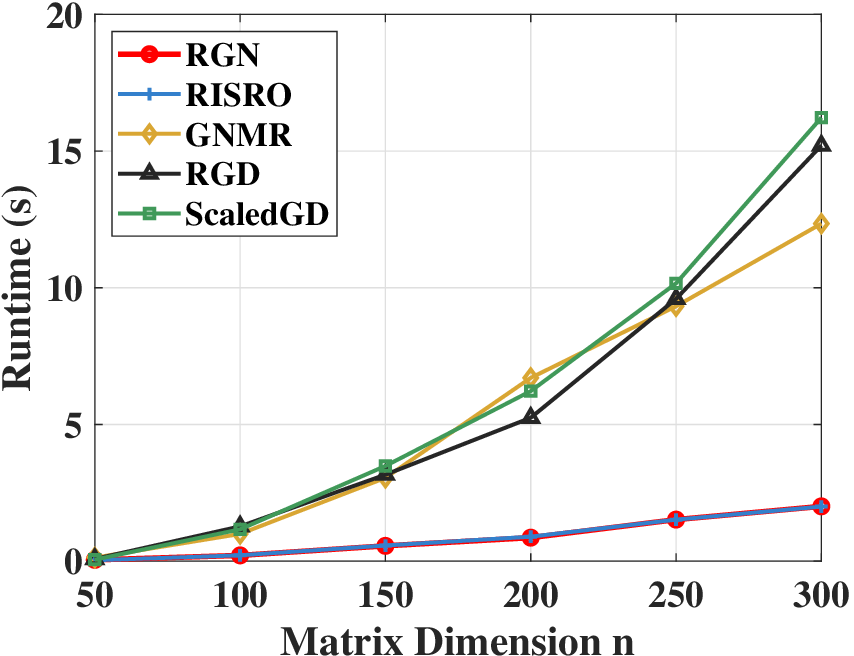}
	\end{minipage}
	\hfill
	\begin{minipage}{0.48\textwidth}
		\centering
		\includegraphics[height=4.5cm]{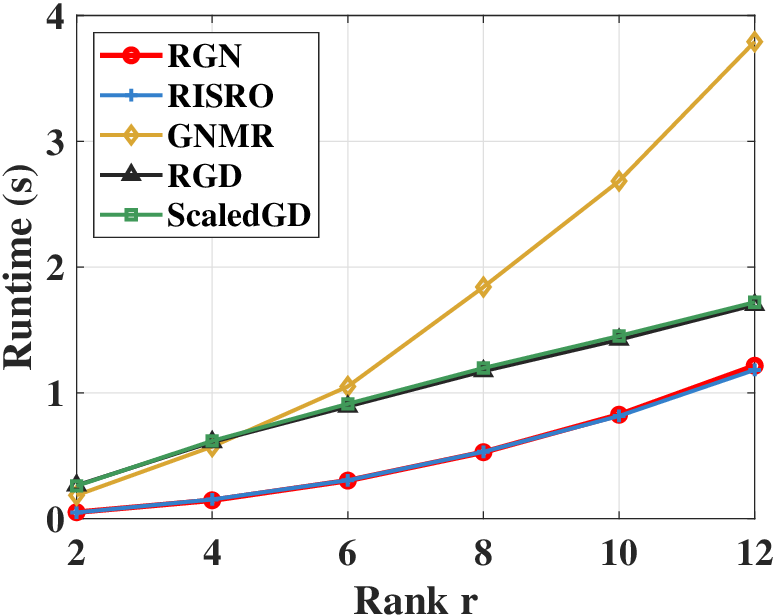}
	\end{minipage}
	\caption{Runtime comparisons.}
	\label{fig:runtime_scaling}
\end{figure}

We first illustrate the effect of the matrix dimension by setting $n_1=n_2=n$, fixing $r=3$, and varying $n$. We then illustrate the effect of the rank by setting $n_1=n_2=100$ and varying $r$. As shown in Figure \ref{fig:runtime_scaling}, the runtime of RGN scales similarly to that of RISRO and remains consistently lower than those of GNMR, RGD, and ScaledGD as either $n$ or $r$ increases. These results demonstrate the computational efficiency and scalability of RGN for low-rank matrix sensing.

\section{Conclusion} \label{sec:conclusion}
In this paper, we have shown that under Gaussian sensing, RGN can recover the ground-truth rank-$r$ matrix with quadratic convergence provided the sample size satisfies $m = \cO ( r(n_1+n_2) \kappa^2)$. This establishes, for the first time, that RGN achieves both rank-optimal sample complexity and quadratic convergence. The key ingredient in our proof is an exact finite Neumann expansion combined with the virtual-sequence decoupling technique, which enables a sharp spectral-norm analysis of the RGN update.

An important direction for future work is to improve the rank-optimal sample complexity to the optimal sample complexity by removing the dependence on the condition number $\kappa$. The phase transition experiment in Figure \ref{fig:phase_transition}(f) suggests that successful recovery is largely insensitive to $\kappa$, while the $\kappa^2$ dependence in our theoretical guarantee arises solely from the analysis of spectral initialization in Lemma \ref{lem:error_init}. This suggests that such dependence may be an artifact of the analysis rather than an intrinsic requirement of RGN.

\appendix

\section{Preliminary Lemmas}
\begin{lemma}[{\cite[Lemma 4.1]{wei2016guarantees}\cite[Lemma 12]{cai2025fast}}] \label{lem:wei}
	Let $\cP_T$ be the orthogonal projection onto the tangent space $T$ of the rank-$r$ matrix manifold at $\m{X}$. Let $\m{Z}$ be another rank-$r$ matrix. Then,
	\equs{
		\frobn{ \sbra{\cI - \cP_T} \m{Z} } \le \frac{\frobn{\m{X} - \m{Z}}^2}{\sigma_r \sbra{ \m{Z} }}, \qquad \twon{ \sbra{\cI - \cP_T}  \m{Z} } \le \frac{\twon{\m{X} - \m{Z} }^2}{\sigma_r \sbra{ \m{Z} }}.
	}
\end{lemma}

\begin{lemma}[{\cite[Lemma 11]{cai2025fast}}] \label{lem:twoProj}
	Let $\cP_{T_i}$ be the orthogonal projection onto the tangent space $T_i$ of the rank-$r$ matrix manifold at $\m{X}_i$ for $i=1,2$. Then,
	\equs{
		\frobn{ (\cP_{T_1} - \cP_{T_2}) \m{Z} } \le \frac{2\sqrt{2} \frobn{\m{X}_1 - \m{X}_2 }}{\sigma_r (\m{X}_1)} \twon{\m{Z}}.
	}
\end{lemma}

\begin{lemma}[{\cite[Lemma 10]{cai2025fast}}] \label{lem:TrTr}
	Let $\m{Z}_1,\m{Z}_2\in\bR^{n_1\times n_2}$ and $\eta>0$.
	Suppose that
	\equs{
		\begin{split}
			\sigma_r(\m{Z}_1)
			& > \sigma_{r+1}(\m{Z}_1), \\
			\sigma_r(\m{Z}_1)-\sigma_{r+1}(\m{Z}_1)
			& \ge \eta \, \sigma_{r+1}(\m{Z}_1), \\
			\twon{\m{Z}_1-\m{Z}_2}
			& \leq (1-1/\sqrt{2}) \sbra{\sigma_r(\m{Z}_1)-\sigma_{r+1}(\m{Z}_1)}.
		\end{split}
	}
	Then,
	\equs{
		\frobn{\cT_r(\m{Z}_1)-\cT_r(\m{Z}_2)}
		\leq (12\eta^{-1} +20) \frobn{\m{Z}_1-\m{Z}_2}.
	}
\end{lemma}

\begin{lemma}[{\cite[Lemma 13]{cai2025fast}}] \label{lem:decoupling}
	Define $\tilde{\cP}\m{Z}:=\langle \m{x}\m{y}^T,\m{Z}\rangle\m{x}\m{y}^T$ and $\tilde{\cP}^{\perp}:=\cI-\tilde{\cP}$.
	For any $\m{Z}\in\bR^{n_1\times n_2}$, we have
	\equs{
		(\cA^*\cA-\tilde{\cA}^*\tilde{\cA})\m{Z}
		=
		\sbra{\cA^*\cA-\cI}\tilde{\cP}\m{Z}
		+
		\langle \m{x}\m{y}^T,\cA^*\cA\tilde{\cP}^{\perp}\m{Z}\rangle\m{x}\m{y}^T.
	}
\end{lemma}

\begin{lemma}[{\cite[Lemma 3]{cai2025fast}}] \label{lem:directional-decoupling}
	Fix $\m{x}\in\bS^{n_1-1}$ and $\m{y}\in\bS^{n_2-1}$. Define $\tilde{\cP}\m{Z}:=\langle \m{x}\m{y}^T,\m{Z}\rangle\m{x}\m{y}^T$ and $\tilde{\cP}^{\perp}:=\cI-\tilde{\cP}$.
	Let $\m{Z}$ is independent of
	$\{\tilde{\cP} \m{A}_i\}_{i=1}^m$.
	Then, with probability at least $1-2\exp\sbra{-8(n_1+n_2)}$, we have
	\equs{
		\abs{ \langle \m{x}\m{y}^T,\cA^*\cA \tilde{\cP}^{\perp} \m{Z}\rangle }
		\le
		4\sqrt{\frac{n_1+n_2}{m}}\twon{ \cA\tilde{\cP}^{\perp}\m{Z} }.
	}
\end{lemma}

The following lemma is a direct result of \cite[Proof of Lemma 5]{cai2025fast}.

\begin{lemma} \label{lem:spectral-control}
	Let $\cN$ be a $\frac{1}{4}$-net of $\bS^{n_1-1}\times\bS^{n_2-1}$.
	Let $\m{M}$ be a matrix derived from the RGN iterates of with $\rank(\m{M}) \le 2r$.
	For each $(\m{x},\m{y})\in\cN$, let
	$\tilde{\m{M}}$ be its counterpart derived from the virtual sequence and satisfying $\rank(\tilde{\m{M}})\le 2r$.
	Suppose that the RIP constant $\delta_{6r} \le 1$. Then, with probability at least $1-2\exp(-5(n_1+n_2))$,
	\equs{
			\twon{ \sbra{\cA^*\cA-\cI}\m{M} }
			\le
			\alpha_1\twon{\m{M}}
			+
			\alpha_2\sup_{(\m{x},\m{y})\in\cN}\frobn{ \m{M} - \tilde{\m{M}} },
		}
	where $\alpha_1:=16\sqrt{2r(n_1+n_2)/m}+2\delta_{6r}$ and $\alpha_2:=4\delta_{6r}+16\sqrt{(n_1+n_2)/m}$. 
\end{lemma}

\section{Technical Lemmas}

\begin{lemma} \label{lem:RIP_tildeA}
	Let $\delta_r$ and $\tilde{\delta}_r$ denote the rank-$r$ RIP constants of $\cA$ in \eqref{eq:sensing} and $\tilde{\cA}$ in \eqref{eq:virtual-sensing-operator}, respectively.
	Then,
	\equs{
		\tilde{\delta}_r \le \delta_{r+1}.
	}
\end{lemma}

\begin{proof}
	Let $\m{Z}\in\bR^{n_1\times n_2}$ satisfy $\rank(\m{Z})\le r$.
	By the definition of $\tilde{\cA}$,
	\equs{
			\twon{\tilde{\cA}\m{Z}}^2
			=
			\twon{\cA(\tilde{\cP}^{\perp}\m{Z})}^2
			+
			\frobn{\tilde{\cP}\m{Z}}^2.
	}
	By the orthogonality of $\tilde{\cP}$ and $\tilde{\cP}^{\perp}$,
	\equs{
			\frobn{\m{Z}}^2
			=
			\frobn{\tilde{\cP}^{\perp}\m{Z}}^2
			+
			\frobn{\tilde{\cP}\m{Z}}^2.
	}
	Consequently,
	\equs{
		\left | \twon{ \tilde{\cA} \m{Z} }^2 - \frobn{ \m{Z} }^2 \right | 
		= \left | \twon{ \cA( \tilde{\cP}^\perp \m{Z} ) }^2 - \frobn{ \tilde{\cP}^\perp \m{Z} }^2 \right | 
		\le \delta_{r +1} \frobn{ \tilde{\cP}^\perp \m{Z} }^2 
		\le \delta_{r +1} \frobn{ \m{Z} }^2, 
	}
	where the first inequality uses the RIP of $\cA$ and $\rank( \tilde{\cP}^\perp \m{Z} ) \le r +1$. 
	Thus, $\delta_{r+1}$ is an admissible rank-$r$ RIP constant for $\tilde{\cA}$, and hence $\tilde{\delta}_r\le\delta_{r+1}$.
\end{proof}

By combining Lemma \ref{lem:RIP_tildeA}, \cite[equation (9)]{cai2025fast}, and Lemma \ref{lem:tangent_RIP}, we obtain the following lemma.
\begin{lemma} \label{lem:PT_I_AA_PT}
	Let $\delta_{2r}$ be the rank-$2r$ RIP constant of $\cA$, and let $\cP_T$ be the orthogonal projection onto the tangent space $T$ at a rank-$r$ matrix.
	Then,
	\equs{
		\begin{split}
			\twon{ \cP_T (\cI - \cA^* \cA) \cP_T } & \le \delta_{2r}. \\
			\twon{ \cP_T (\cI - \tilde{\cA}^* \tilde{\cA}) \cP_T } & \le \delta_{2r+1}.
		\end{split}
	}
	Further, if $\delta_{2r+1} < 1$, then
	\equs{
		\frac{1}{1+\delta_{2r+1}} \frobn{\m{Z}}
		\le \frobn{ [ \cP_T \tilde{\cA}^* \tilde{\cA} \cP_T ]^{-1} \m{Z} }
		\le \frac{1}{1-\delta_{2r+1}} \frobn{\m{Z}}.
	}
\end{lemma}

The following Lemma \ref{lem:RIP-inner-product} extends \cite[Lemma 4.3]{wei2016guarantees} by removing the orthogonality assumption $\left\langle \m{X}_1,\m{X}_2 \right\rangle=0$.
\begin{lemma} \label{lem:RIP-inner-product}
	Let $\m{X}_1,\m{X}_2 \in \bR^{n_1\times n_2}$. Suppose
	$\rank \sbra{\m{X}_1}+\rank \sbra{\m{X}_2} \le \min (n_1,n_2)$.
	Then,
	\equs{
		\abs{
			\left\langle
			\m{X}_1,
			\sbra{\cA^*\cA-\cI}\m{X}_2
			\right\rangle
		}
		\le
		\delta_{\rank \sbra{\m{X}_1}+\rank \sbra{\m{X}_2}}
		\frobn{\m{X}_1}
		\frobn{\m{X}_2}.
	}
\end{lemma}

\begin{proof}
	Set $q:=\rank \sbra{\m{X}_1}+\rank \sbra{\m{X}_2}$, $\alpha:=\frobn{\m{X}_1}$, and $\beta:=\frobn{\m{X}_2}$.
	If $\alpha\beta=0$, then the result holds immediately.
	When $\alpha\beta\neq 0$, we have
	\equs{
		\begin{split}
			\left\langle
			\cA(\m{X}_1),
			\cA(\m{X}_2)
			\right\rangle
			& =
			\frac{\alpha\beta}{4}
			\mbra{
				\twon{
					\cA\sbra{
						\alpha^{-1}\m{X}_1
						+
						\beta^{-1}\m{X}_2
					}
				}^2
				-
				\twon{
					\cA\sbra{
						\alpha^{-1}\m{X}_1
						-
						\beta^{-1}\m{X}_2
					}
				}^2
			} \\
			& \le
			\frac{\alpha\beta}{4}
			\mbra{
				(1+\delta_q)
				\frobn{
					\alpha^{-1}\m{X}_1
					+
					\beta^{-1}\m{X}_2
				}^2
				-
				(1-\delta_q)
				\frobn{
					\alpha^{-1}\m{X}_1
					-
					\beta^{-1}\m{X}_2
				}^2
			} \\
			& =
			\left\langle \m{X}_1,\m{X}_2 \right\rangle
			+
			\frac{\alpha\beta\delta_q}{2}
			\mbra{
				\frobn{\alpha^{-1}\m{X}_1}^2
				+
				\frobn{\beta^{-1}\m{X}_2}^2
			} \\
			& =
			\left\langle \m{X}_1,\m{X}_2 \right\rangle
			+
			\delta_q
			\frobn{\m{X}_1}
			\frobn{\m{X}_2},
		\end{split}
	}
	where the inequality follows from the RIP of $\cA$ and $\rank \sbra{ \alpha^{-1}\m{X}_1 \pm \beta^{-1}\m{X}_2} \le q$.
	Applying the same argument with $\m{X}_2$ replaced by $-\m{X}_2$ yields
	\equs{
		\left\langle
		\cA(\m{X}_1),
		\cA(\m{X}_2)
		\right\rangle
		\ge
		\left\langle \m{X}_1,\m{X}_2 \right\rangle
		-
		\delta_q
		\frobn{\m{X}_1}
		\frobn{\m{X}_2}.
	}
	Therefore,
	\equs{
		\abs{
			\left\langle
			\m{X}_1,
			\sbra{\cA^*\cA-\cI}\m{X}_2
			\right\rangle
		}
		=
		\abs{
			\left\langle
			\cA(\m{X}_1),
			\cA(\m{X}_2)
			\right\rangle
			-
			\left\langle \m{X}_1,\m{X}_2 \right\rangle
		}  \le
		\delta_q
		\frobn{\m{X}_1}
		\frobn{\m{X}_2}.
	}
\end{proof}

\begin{lemma} \label{lem:proj-RIP}
	Let $\rank \sbra{\m{Z}} \le r_1$, and let $\cP_T$ be the orthogonal projection onto the tangent space $T$ of the rank-$r$ matrix manifold. Then, 
	\equs{
		\begin{split}
			\twon{ \sbra{\cA^*\cA-\cI}\m{Z} }
			& \le \delta_{r_1+1} \frobn{\m{Z}}, \\
			\frobn{ \cP_T \sbra{\cA^*\cA-\cI}\m{Z} }
			& \le \delta_{r_1+2r} \frobn{\m{Z}}.
		\end{split}
	}
\end{lemma}

\begin{proof}
	Applying Lemma \ref{lem:RIP-inner-product}, we have
	\equs{
		\twon{ \sbra{\cA^*\cA-\cI}\m{Z} }
		= \sup_{\twon{\m{x}}=\twon{\m{y}}=1}
		\abs{
			\langle
			\m{x}\m{y}^T,
			\sbra{\cA^*\cA-\cI}\m{Z}
			\rangle
		} 
		\le \delta_{r_1+1}\frobn{\m{Z}},
	}
	and
	\equs{
		\begin{split}
			\frobn{
				\cP_T\sbra{\cA^*\cA-\cI}\m{Z}
			}
			& =
			\sup_{\frobn{\m{Y}}=1}
			\abs{
				\left\langle
				\m{Y},
				\cP_T\sbra{\cA^*\cA-\cI}\m{Z}
				\right\rangle
			} \\
			& =
			\sup_{\frobn{\m{Y}}=1}
			\abs{
				\left\langle
				\cP_T \m{Y},
				\sbra{\cA^*\cA-\cI}\m{Z}
				\right\rangle
			} \\
			& \le
			\delta_{r_1+2r}\frobn{\m{Z}}.
		\end{split}
	}
\end{proof}

\section{Proof of Lemma \ref{lem:tangent_RIP}} \label{append:tangent_RIP}
	The operator $\cP_T \cA^*\cA \cP_T$ is self-adjoint on $T$. To show the first line, it is equivalent to show the spectrum of $\cP_T \cA^*\cA \cP_T$ lies in $[1-\delta_{2r}, 1+\delta_{2r}]$. In particular, its eigenvalues admit the variational characterizations
	\equs{
		\lambda_{\max}
		:= \sup_{\m{Z}\in T, \frobn{\m{Z}}=1}
		\left\langle \cP_T \cA^*\cA \cP_T \m{Z},\m{Z} \right\rangle, \qquad
		\lambda_{\min}
		:= \inf_{\m{Z}\in T, \frobn{\m{Z}}=1}
		\left\langle \cP_T \cA^*\cA \cP_T \m{Z},\m{Z} \right\rangle.
	}
	For every $\m{Z}\in T$, the identities $\cP_T\m{Z}=\m{Z}$ and $\cP_T^*=\cP_T$ give
	\equ{
		\left\langle \cP_T \cA^*\cA \cP_T \m{Z},\m{Z} \right\rangle
		= \left\langle \cA^*\cA \m{Z},\m{Z} \right\rangle 
		= \left\langle \cA \m{Z},\cA \m{Z} \right\rangle 
		= \twon{ \cA \m{Z} }^2.
		\label{eq:tangent-quadratic}
	}
	Every matrix in $T$ has rank at most $2r$.
	Therefore, the rank-$2r$ RIP and \eqref{eq:tangent-quadratic} imply
	\equs{
		(1-\delta_{2r}) \frobn{\m{Z}}^2
		\le
		\left\langle \cP_T \cA^*\cA \cP_T \m{Z},\m{Z} \right\rangle
		\le
		(1+\delta_{2r}) \frobn{\m{Z}}^2.
	}
	This yields
	\equs{
		1-\delta_{2r}
		\le
		\lambda_{\min}
		\le
		\lambda_{\max}
		\le
		1+\delta_{2r}.
	}
	Since $\delta_{2r}<1$ and $\lambda_{\min} > 0$, the $\cP_T \cA^*\cA \cP_T$ is positive definite and hence invertible on $T$.
	For any $\m{X}\in T$, let $\m{Z} = [\cP_T \cA^*\cA \cP_T]^{-1} \m{X} \in T$. Substituting $\m{Z}$ into the first line gives the second line.

\section{Proof of Lemma \ref{lem:error_contraction_F}} \label{append:error_contraction_F}
We have
\equs{
	\begin{split}
		& \frobn{ \m{X}_{t+1} - \m{X}_\star } \\
		& \le
		\frobn{\m{X}_{t+1} - \m{W}_{t+1}}
		+ \frobn{\m{W}_{t+1} - \m{X}_{\star}}
		\\
		& \le 2 \frobn{\m{W}_{t+1} - \m{X}_{\star}} \\
		& \le 2 \frobn{ (\cP_{T_t} - \cI) \m{X}_{\star} } + 2 \frobn{ \mbra{\cP_{T_t} \cA^*\cA \cP_{T_t}}^{-1}
			\cP_{T_t} \sbra{\cA^*\cA - \cI}
			(\cI - \cP_{T_t}) \m{X}_{\star} } \\
		& \le \sbra{ 2 + \frac{2\delta_{3r}}{1-\delta_{2r}} } \frobn{ (\cP_{T_t} - \cI) \m{X}_{\star} } \\
		& \le \frac{4}{\sigma_r(\m{X}_\star)} \frobn{ \m{X}_t - \m{X}_\star }^2,
	\end{split} 
}
where the second inequality holds since $\rank \sbra{\m{X}_{\star}} = r$ and $\m{X}_{t+1}$ is a best rank-$r$ approximation of $\m{W}_{t+1}$, the third follows from Lemma \ref{lem:tangent_RIP} and Lemma \ref{lem:proj-RIP}, and the last follows from $\delta_{2r} \le \delta_{3r} \le 1/2$ and Lemma \ref{lem:wei}.

By induction and \eqref{eq:XT_Xstar}, we have
\equs{
	\begin{split}
		\frobn{ \m{X}_t - \m{X}_\star } 
		& \le \frac{\sigma_r(\m{X}_\star)}{4} \sbra{ \frac{4}{\sigma_r(\m{X}_\star)} \frobn{ \m{X}_T - \m{X}_\star } }^{2^{t-T}} \\
		& \le 2^{-2^{t-T}} \sigma_r(\m{X}_\star).
	\end{split} 
}

\section{Proof of Lemma \ref{lem:ZminusZtilde}} \label{append:ZminusZtilde}
	Fix $(\m{x},\m{y})\in\cN$. We have
	\equ{
		\begin{split}
			\m{Z}-\tilde{\m{Z}}
			&=
			\cP_{T_t}
			\sbra{\cA^*\cA-\cI}
			(\m{M}-\tilde{\m{M}}) 
			+
			\cP_{T_t}
			(\cA^*\cA-\tilde{\cA}^*\tilde{\cA})
			\tilde{\m{M}} \\
			&\quad +
			(\cP_{T_t}-\cP_{\tilde{T}_t})
			(\tilde{\cA}^*\tilde{\cA}-\cI)
			\tilde{\m{M}} \\
			&:=
			\m{J}_1+\m{J}_2+\m{J}_3.
		\end{split}
		\label{eq:Z-diff}
	}
	
	For $\m{J}_1$, it follows from Lemma \ref{lem:proj-RIP} and $\rank (\m{M}-\tilde{\m{M}}) \le 4r$ that
	\equ{
		\frobn{\m{J}_1}
		\le
		\delta_{6r}
		\frobn{ \m{M}-\tilde{\m{M}} }.
		\label{eq:J1bound}
	}
	
	For $\m{J}_2$, we have
	\equ{
		\begin{split}
			\frobn{\m{J}_2}
			&\le
			\frobn{
				\cP_{T_t}
				\sbra{\cA^*\cA-\cI}
				\tilde{\cP}\tilde{\m{M}}
			} +
			\abs{
				\langle
				\m{x}\m{y}^T,
				\cA^*\cA
				\tilde{\cP}^{\perp}\tilde{\m{M}}
				\rangle
			} \cdot
			\frobn{
				\cP_{T_t}(\m{x}\m{y}^T)
			} \\
			&\le
			\delta_{1+2r}
			\frobn{\tilde{\cP}\tilde{\m{M}}}
			+
			4\sqrt{\frac{n_1+n_2}{m}}
			\twon{
				\cA\tilde{\cP}^{\perp}\tilde{\m{M}}
			} \\
			&\le
			\sbra{
				\delta_{1+2r}
				+
				8\sqrt{\frac{r(n_1+n_2)}{m}}
			}
			\twon{\tilde{\m{M}}} \\
			& \le
			\frac{\alpha_1}{2}
			( \twon{\m{M}} + \frobn{ \m{M}-\tilde{\m{M}} } ).
		\end{split}
		\label{eq:J2bound}
	}
	Here the first inequality uses Lemma \ref{lem:decoupling}. Since $\tilde{\m{M}}$ is generated from the virtual sequence, it is independent of $\{\tilde{\cP} \m{A}_i\}_{i=1}^m$, allowing us to invoke Lemma \ref{lem:directional-decoupling} to bound $\abs{ \langle \m{x}\m{y}^T, \cA^*\cA \tilde{\cP}^{\perp}\tilde{\m{M}} \rangle }$. Hence, the second uses Lemma \ref{lem:proj-RIP} with $\rank(\tilde{\cP}\tilde{\m{M}}) = 1$, Lemma \ref{lem:directional-decoupling}, and the fact $\frobn{\cP_{T_t}(\m{x}\m{y}^T)}\le 1$. 
	The third inequality in \eqref{eq:J2bound} uses 
	\equ{
		\begin{split}
			\frobn{\tilde{\cP}\tilde{\m{M}}} 
			& = \abs{ \langle \m{x}\m{y}^T,\tilde{\m{M}} \rangle } \cdot \frobn{ \m{x}\m{y}^T }
			\le \twon{\tilde{\m{M}}}, \\
			\twon{\cA \tilde{\cP}^{\perp}\tilde{\m{M}} } 
			& \le \sqrt{1+\delta_{2r+1}} \frobn{ \tilde{\cP}^{\perp}\tilde{\m{M}} } \le \sqrt{1+\delta_{2r+1}} \frobn{ \tilde{\m{M}} } \le 2\sqrt{r} \twon{ \tilde{\m{M}} },
		\end{split} \label{eq:PM_APM}
	}
	where we use $\delta_{2r+1} \le 1$, and the last uses the definition of $\alpha_1$ and the triangle inequality.
	
	For $\m{J}_3$, Lemma \ref{lem:twoProj} gives
	\equ{
		\frobn{\m{J}_3}
		\le
		\frac{
			2\sqrt{2}
			\frobn{\m{X}_t-\tilde{\m{X}}_t}
		}{
			\sigma_r(\m{X}_t)
		}
		\twon{
			(\tilde{\cA}^*\tilde{\cA}-\cI)
			\tilde{\m{M}}
		}.
		\label{eq:J3}
	}
	We have the decomposition
	\equ{
		\begin{split}
			\twon{
				(\tilde{\cA}^*\tilde{\cA}-\cI)
				\tilde{\m{M}}
			} 
			&\le
			\twon{
				\sbra{\cA^*\cA-\cI}\m{M}
			}
			+
			\twon{
				\sbra{\cA^*\cA-\cI}
				(\tilde{\m{M}}-\m{M})
			} \\
			&\quad+
			\twon{
				(\cA^*\cA-\tilde{\cA}^*\tilde{\cA})
				\tilde{\m{M}}
			} \\
			&:=
			J_{31}+J_{32}+J_{33}.
		\end{split}
		\label{eq:Z0-vs}
	}
	For $J_{31}$, it follows from Lemma \ref{lem:spectral-control} that
	\equs{
		J_{31}
		\le
		\alpha_1 \twon{\m{M}}
		+
		\alpha_2 \sup_{(\m{x},\m{y})\in\cN}
		\frobn{\m{M}-\tilde{\m{M}}}.
	}
	For $J_{32}$, it follows from Lemma
	\ref{lem:proj-RIP} that
	\equs{
		\begin{split}
			J_{32}
			\le
			\delta_{4r+1}
			\frobn{ \m{M} - \tilde{\m{M}}}.
		\end{split}
	}
	
	\noindent For $J_{33}$, we have
	\equs{
		\begin{split}
			J_{33}
			&\le
			\twon{
				\sbra{\cA^*\cA-\cI}
				\tilde{\cP}\tilde{\m{M}}
			} +
			\abs{
				\langle
				\m{x}\m{y}^T,
				\cA^*\cA
				\tilde{\cP}^{\perp}\tilde{\m{M}}
				\rangle
			} \cdot
			\twon{\m{x}\m{y}^T} \\
			&\le
			\delta_{2}\frobn{\tilde{\cP}\tilde{\m{M}}}
			+
			4\sqrt{\frac{n_1+n_2}{m}}
			\twon{
				\cA\tilde{\cP}^{\perp}\tilde{\m{M}}
			} \\
			&\le
			\sbra{
				\delta_2
				+
				8 \sqrt{\frac{r(n_1+n_2)}{m}}
			}
			\twon{\tilde{\m{M}}} \\
			& \le \frac{\alpha_1}{2} ( \twon{\m{M}} + \frobn{ \m{M}-\tilde{\m{M}} } ).
		\end{split}
	}
	Here the first inequality uses Lemma \ref{lem:decoupling}. The second uses Lemma \ref{lem:proj-RIP} and Lemma \ref{lem:directional-decoupling}. The third uses \eqref{eq:PM_APM}.
	
	Substituting the bounds of $J_{31}$, $J_{32}$, and
	$J_{33}$ into \eqref{eq:Z0-vs} and then into \eqref{eq:J3} gives
	\equ{
		\begin{split}
			\frobn{\m{J}_3}
			&\le
			\frac{
				2\sqrt{2}
				\frobn{\m{X}_t-\tilde{\m{X}}_t}
			}{
				\sigma_r(\m{X}_t)
			}
			\Bigg(
			\frac{3\alpha_1}{2}\twon{\m{M}} +
			3 \alpha_1
			\sup_{(\m{x},\m{y})\in\cN}
			\frobn{\m{M}-\tilde{\m{M}}}
			\Bigg) \\
			& \le \frac{9 c_0 \alpha_1}{1-c_0} \Bigg( \twon{\m{M}} +
			\sup_{(\m{x},\m{y})\in\cN}
			\frobn{\m{M}-\tilde{\m{M}}} \Bigg).
		\end{split}
		\label{eq:J3bound}
	}
	Here the first inequality follows from $\alpha_2 \le 2 \alpha_1$ and $\delta_{4r+1} \le \delta_{6r} \le \alpha_1/2$. The second uses $\frobn{\m{X}_t-\tilde{\m{X}}_t} \le E_t \le c_0 \sigma_r (\m{X}_\star)$ and 
	\equ{
		\sigma_r(\m{X}_t) 
		\ge \sigma_r(\m{X}_\star) - \twon{ \m{X}_\star - \m{X}_t } \ge \sigma_r(\m{X}_\star) - E_t \ge (1-c_0) \sigma_r(\m{X}_\star), \label{eq:sigma_r_Xt}
	}
	where we use Weyl's inequality.
	
	Substituting \eqref{eq:J1bound}, \eqref{eq:J2bound} and \eqref{eq:J3bound} into \eqref{eq:Z-diff}, and then taking the supremum over $(\m{x},\m{y})\in\cN$, yields the desired estimate.
	
	It remains to track the failure probability.
	Lemma \ref{lem:directional-decoupling} is applied for each fixed $(\m{x},\m{y})\in\cN$, whereas Lemma \ref{lem:spectral-control} holds uniformly over $\cN$.
	Therefore, a union bound and $\abs{\cN}\le 12^{n_1+n_2}$ give
	\equs{
		\begin{split}
			2 \abs{\cN}\exp\sbra{-8(n_1+n_2)}
			+
			2\exp\sbra{-5(n_1+n_2)}
			\le
			4 \exp\sbra{-5(n_1+n_2)}.
		\end{split}
	}
	Hence, the desired estimate holds with probability at least $1 - 4 \exp\sbra{-5(n_1+n_2)}$.

\section{Proof of Lemma \ref{lem:error_contraction}} \label{append:error_contraction}

Recalling the definition of $E_t$ in \eqref{eq:Et}, we separately bound
$\twon{\m{X}_{t+1}-\m{X}_{\star}}$ and
$\sup_{(\m{x},\m{y})\in\cN}\frobn{\m{X}_{t+1}-\tilde{\m{X}}_{t+1}}$.

Using \eqref{eq:neumann} and \eqref{eq:Z0Zk}, we rewrite $\m{W}_{t+1}$ in \eqref{eq:W} as
\equ{
	\begin{split}
		\m{W}_{t+1} 
		= \cP_{T_t} \m{X}_{\star} 
		+ \sum_{k=0}^{K-1} \m{Z}_k 
		+ \mbra{\cP_{T_t} \sbra{\cI - \cA^*\cA} \cP_{T_t}}^K
		\mbra{\cP_{T_t} \cA^*\cA \cP_{T_t}}^{-1} \m{Z}_0.
	\end{split}
	\label{eq:Wnew}
}
Similarly, using \eqref{eq:neumann} and \eqref{eq:Z0Zk_tilde}, we rewrite $\tilde{\m{W}}_{t+1} $ as
\equ{
	\begin{split}
		\tilde{\m{W}}_{t+1} 
		= \cP_{\tilde{T}_t} \m{X}_{\star} 
		+ \sum_{k=0}^{K-1} \tilde{\m{Z}}_k 
		+ [ \cP_{\tilde{T}_t} ( \cI - \tilde{\cA}^*\tilde{\cA} ) \cP_{\tilde{T}_t} ]^K
		[ \cP_{\tilde{T}_t} \tilde{\cA}^*\tilde{\cA} \cP_{\tilde{T}_t} ]^{-1} \tilde{\m{Z}}_0.
	\end{split}
	\label{eq:Wtildenew}
}

\subsection{
	Bounding
	\texorpdfstring{
		$
		\displaystyle
		\twon{\m{X}_{t+1} - \m{X}_{\star}}
		$
	}{
		x
	}
}
\noindent

It follows from \eqref{eq:Xt_Xstar_2} that
\equ{
		\twon{\m{X}_{t+1} - \m{X}_{\star}}
		\le 2 \twon{\m{W}_{t+1} - \m{X}_{\star}}.
	\label{eq:truncation}
}
It follows from \eqref{eq:Wnew} that
\equ{
	\begin{split}
		\twon{\m{W}_{t+1} - \m{X}_{\star}}
		& \le \twon{(\cP_{T_t} - \cI) \m{X}_{\star}}
		+ \sum_{k=0}^{K-1} \twon{\m{Z}_k}
		\\
		& \quad + \twon{\mbra{\cP_{T_t} \sbra{\cI - \cA^*\cA} \cP_{T_t}}^K
			\mbra{\cP_{T_t} \cA^*\cA \cP_{T_t}}^{-1}
			\m{Z}_0 } \\
		& := I_1 + I_2 + I_3.
	\end{split}
	\label{eq:error-terms}
}
\subsubsection{Bounding $I_1$}
It follows from Lemma \ref{lem:wei} that gives
\equ{
	I_1 \le \frac{\twon{\m{X}_t - \m{X}_{\star}}^2}{\sigma_r (\m{X}_{\star})} \le \frac{ E_t^2 }{\sigma_r (\m{X}_{\star})}.
	\label{eq:I1}
}

\subsubsection{Bounding $I_2$}
We bound $\twon{\m{Z}_k}$ by induction.

\medskip
\noindent\textit{(i) Bounding $\twon{\m{Z}_0}$.}
We have
\equ{
	\begin{split}
		\twon{\m{Z}_0} 
		& = \twon{ \cP_{T_t} \sbra{\cA^*\cA - \cI} (\cI - \cP_{T_t}) \m{X}_{\star} } \\
		& \le 2 \twon{ \sbra{\cA^*\cA - \cI} (\cI - \cP_{T_t}) \m{X}_{\star} } \\
		& \le 2 \alpha_1 \twon{ (\cI - \cP_{T_t}) \m{X}_{\star} } + 2 \alpha_2 \sup_{(\m{x},\m{y}) \in \cN} \frobn{ (\cP_{T_t} - \cP_{\tilde{T}_t}) \m{X}_{\star} }.
	\end{split} \label{eq:Z0twonorm}
}
Here the first inequality uses $\twon{\cP_{T_t} \m{Y} } = \twon{ \m{U}_t\m{U}_t^H \m{Y} \sbra{ \m{I} - \m{V}_t\m{V}^H_t } + \m{Y} \m{V}_t\m{V}^H_t } \le 2 \twon{\m{Y}}$ for any $\m{Y}$. The second uses Lemma \ref{lem:spectral-control} by letting $\m{M} = (\cI-\cP_{T_t})\m{X}_\star$ and $\tilde{\m{M}} = (\cI-\cP_{\tilde{T}_t})\m{X}_\star$.

In \eqref{eq:Z0twonorm}, the first term can be bounded using Lemma \ref{lem:wei}, so we focus on the second term. Fix $(\m{x},\m{y}) \in \cN$.
Then,
\equ{
	\begin{split}
		\frobn{ (\cP_{T_t} - \cP_{\tilde{T}_t}) \m{X}_{\star} }
		& = \frobn{ (\cP_{T_t} - \cP_{\tilde{T}_t}) (\m{X}_\star - \m{X}_t) + (\cP_{T_t} - \cP_{\tilde{T}_t}) \m{X}_t } \\
		& \le \frobn{ (\cP_{T_t} - \cP_{\tilde{T}_t} ) (\m{X}_\star - \m{X}_t) } + \frobn{ ( \cI - \cP_{\tilde{T}_t} ) \m{X}_t } \\
		& \le \frac{2\sqrt{2} \frobn{ \m{X}_t - \tilde{\m{X}}_t }}{\sigma_r(\m{X}_t)} \twon{ \m{X}_\star - \m{X}_t } + \frac{1}{\sigma_r(\m{X}_t)} \frobn{ \m{X}_t - \tilde{\m{X}}_t }^2 \\
		& \le \frac{2\sqrt{2}+1}{1-c_0} \frac{E_t^2}{\sigma_r(\m{X}_\star)} \\
		& \le \frac{5 E_t^2}{\sigma_r(\m{X}_\star)}. \label{eq:PTminusPT}
	\end{split}
}
Here the first inequality follows from $\cP_{T_t} \m{X}_t = \m{X}_t$, the second follows from Lemma \ref{lem:twoProj} and Lemma \ref{lem:wei}, the third follows from \eqref{eq:sigma_r_Xt} and the definition of $E_t$, and the last uses $0<c_0<10^{-4}$.

Substituting \eqref{eq:PTminusPT} into \eqref{eq:Z0twonorm} and using Lemma \ref{lem:wei}, we obtain
\equ{
	\twon{\m{Z}_0} 
	\le 2 \alpha_1 \frac{\twon{\m{X}_t - \m{X}_\star}^2}{\sigma_r(\m{X}_\star)} + 10 \alpha_2 \frac{E_t^2}{\sigma_r(\m{X}_\star)} 
	\le  \frac{ 2 E_t^2}{\sigma_r(\m{X}_\star)}.
	\label{eq:Z0}
}
Here the last inequality uses $\alpha_2 \le 2 \alpha_1$ and $\alpha_1 \le 1/12$, which is guaranteed by the assumption $m \ge C r(n_1+n_2)$ for a sufficiently large $C>0$.

\medskip
\noindent\textit{(ii) Bounding $\sup_{(\m{x},\m{y}) \in \cN} \frobn{\m{Z}_0 - \tilde{\m{Z}}_0}$.}
We have
\equ{
	\begin{split}
	& \sup_{(\m{x},\m{y}) \in \cN} \frobn{\m{Z}_0 - \tilde{\m{Z}}_0} \\
	& = \sup_{(\m{x},\m{y}) \in \cN} \frobn{ \cP_{T_t} \sbra{\cA^*\cA - \cI} (\cI - \cP_{T_t}) \m{X}_{\star} - \cP_{\tilde{T}_t} (\tilde{\cA}^*\tilde{\cA} - \cI) (\cI - \cP_{\tilde{T}_t}) \m{X}_{\star} } \\
	& \le \frac{(8 c_0+1)\alpha_1}{1-c_0} \sbra{ \twon{ (\cI-\cP_{T_t})\m{X}_\star } + 
	\sup_{(\m{x},\m{y})\in\cN} \frobn{ (\cP_{T_t} - \cP_{\tilde{T}_t}) \m{X}_{\star} } } \\
	& \le 2 \alpha_1 \sbra{ \frac{\twon{\m{X}_t - \m{X}_{\star}}^2}{\sigma_r (\m{X}_{\star})} + \frac{5 E_t^2}{\sigma_r(\m{X}_\star)} } \\
	& \le  \frac{E_t^2}{\sigma_r(\m{X}_\star)},
	\end{split}
	\label{eq:Z0-F}
}
Here the first inequality follows from Lemma \ref{lem:ZminusZtilde} where we let 
\equs{
	\m{M} = (\cI-\cP_{T_t})\m{X}_\star, \qquad 
	\tilde{\m{M}} = (\cI-\cP_{\tilde{T}_t})\m{X}_\star,
}
the second uses Lemma \ref{lem:wei}, \eqref{eq:PTminusPT} and $0<c_0<10^{-4}$, and the last uses $\alpha_1 \le 1/12$.

\medskip
\noindent\textit{(iii) Bounding $\twon{\m{Z}_{k+1}}$.}
We have
\equ{
	\begin{split}
		\twon{\m{Z}_{k+1}}
		& = \twon{ \cP_{T_t} \sbra{\cI - \cA^*\cA} \m{Z}_k } \\
		& \le 2 \twon{ \sbra{\cA^*\cA - \cI} \m{Z}_k } \\
		& \le 2 \alpha_1 \twon{\m{Z}_k} + 2 \alpha_2 \sup_{(\m{x},\m{y}) \in \cN} \frobn{\m{Z}_k - \tilde{\m{Z}}_k},
	\end{split}
	\label{eq:Zk}
}
where the last inequality follows from Lemma \ref{lem:spectral-control}.

\medskip
\noindent\textit{(iv) Bounding $\sup_{(\m{x},\m{y}) \in \cN} \frobn{\m{Z}_{k+1} - \tilde{\m{Z}}_{k+1}}$. }
We have
\equ{
	\begin{split}
		& \sup_{(\m{x},\m{y}) \in \cN} \frobn{ \m{Z}_{k+1} - \tilde{\m{Z}}_{k+1} } \\
		& = \sup_{(\m{x},\m{y}) \in \cN} \frobn{ \cP_{T_t} \sbra{\cI - \cA^*\cA} \m{Z}_k - \cP_{\tilde{T}_t} (\cI - \tilde{\cA}^*\tilde{\cA}) \tilde{\m{Z}}_k }  \\
		& \le 2 \alpha_1 \sbra{ \twon{ \m{Z}_k } + 
		\sup_{(\m{x},\m{y})\in\cN} \frobn{ \m{Z}_k - \tilde{\m{Z}}_k } } .
	\end{split}
	\label{eq:Zk-F}
}
Here the inequality follows from $0<c_0<10^{-4}$ and Lemma \ref{lem:ZminusZtilde} where we let 
\equs{
	\m{M} = - \m{Z}_k, \qquad 
	\tilde{\m{M}} = - \tilde{\m{Z}}_k.
}

\medskip
\noindent\textit{(v) Final estimate for $I_2$.}
Combining \eqref{eq:Zk}, \eqref{eq:Zk-F}, and $\alpha_1 \le 1/12$ gives the one-step recursion
\equs{
		\twon{\m{Z}_{k+1}} + \sup_{(\m{x},\m{y}) \in \cN} \frobn{\m{Z}_{k+1} - \tilde{\m{Z}}_{k+1}} 
		\le \frac{1}{2} \sbra{ \twon{\m{Z}_k} + \sup_{(\m{x},\m{y}) \in \cN} \frobn{\m{Z}_k - \tilde{\m{Z}}_k} }.
}
By induction, we have
\equ{
	\begin{split}
		\sum_{k=0}^{K-1} \sbra{ \twon{\m{Z}_k} + \sup_{(\m{x},\m{y}) \in \cN} \frobn{ \m{Z}_k-\tilde{\m{Z}}_k } } 
		& \le 2 \sbra{ \twon{\m{Z}_0} + \sup_{(\m{x},\m{y}) \in \cN} \frobn{ \m{Z}_0-\tilde{\m{Z}}_0 } } \\
		& \le \frac{ 6 E_t^2}{\sigma_r(\m{X}_\star)},
	\end{split} \label{eq:finite-sums}
}
where the last inequality uses \eqref{eq:Z0} and \eqref{eq:Z0-F}.
Recall that $I_2 = \sum_{k=0}^{K-1} \twon{\m{Z}_k}$. Consequently,
\equ{
	I_2 \le \frac{6 E_t^2}{\sigma_r(\m{X}_\star)}.
	\label{eq:I2}
}

\subsubsection{Bounding $I_3$}
We have
\equ{
	\begin{split}
		I_3
		& = \twon{ \mbra{\cP_{T_t} \sbra{\cI-\cA^*\cA} \cP_{T_t}}^K \mbra{\cP_{T_t}\cA^*\cA\cP_{T_t}}^{-1} \m{Z}_0 } \\ 
		& \le \frobn{ \mbra{\cP_{T_t} \sbra{\cI-\cA^*\cA} \cP_{T_t}}^K \mbra{\cP_{T_t}\cA^*\cA\cP_{T_t}}^{-1} \m{Z}_0 } \\ 
		& \le \frac{\delta_{2r}^K}{1-\delta_{2r}} \frobn{\m{Z}_0} 
		\le \frac{\sqrt{2r}\,\delta_{2r}^K}{1-\delta_{2r}} \twon{\m{Z}_0}
		\le \frac{6 E_t^2}{\sigma_r(\m{X}_\star)}.
	\end{split}
	\label{eq:I3}
}
Here the second inequality follows from Lemma \ref{lem:PT_I_AA_PT} and Lemma \ref{lem:tangent_RIP}, the third follows from $\rank(\m{Z}_0) \le 2r$, and the last holds since $\delta_{2r} \le 1/2$, \eqref{eq:Z0}, and $\sqrt{r} \delta_{2r}^K \le \sqrt{r} 2^{-\log_2 r} \le r^{-1/2} \le 1$, where we let $K=\max\{1,\lceil\log_2 r\rceil\}$. 

Finally, substituting \eqref{eq:I1}, \eqref{eq:I2}, and \eqref{eq:I3} into \eqref{eq:error-terms} and then into \eqref{eq:truncation} gives
\equ{
	\twon{ \m{X}_{t+1}-\m{X}_\star } 
	\le 2 \twon{\m{W}_{t+1} - \m{X}_{\star}}
	\le \frac{26 E_t^2}{\sigma_r(\m{X}_\star)}.
	\label{eq:X-next}
}

\subsection{
	Bounding
	\texorpdfstring{
		$
		\displaystyle
		\sup_{(\m{x},\m{y})\in\cN}
		\frobn{\m{X}_{t+1}-\tilde{\m{X}}_{t+1}}
		$
	}{
		x
	}
}
\noindent

Note that $\m{X}_{t+1} = \cT_r(\m{W}_{t+1})$ and $\tilde{\m{X}}_{t+1} = \cT_r(\tilde{\m{W}}_{t+1})$. First, to apply Lemma \ref{lem:TrTr}, we need to bound $\frobn{\m{W}_{t+1}-\tilde{\m{W}}_{t+1}} $.
It follows from \eqref{eq:Wnew} and \eqref{eq:Wtildenew} that
\equ{
	\begin{split}
		\frobn{ \m{W}_{t+1}-\tilde{\m{W}}_{t+1} }
		& \le \frobn{ (\cP_{T_t}-\cP_{\tilde{T}_t}) \m{X}_\star } + \sum_{k=0}^{K-1} \frobn{ \m{Z}_k-\tilde{\m{Z}}_k } \\
		& \quad + \frobn{ \mbra{\cP_{T_t} \sbra{\cI-\cA^*\cA} \cP_{T_t}}^K \mbra{\cP_{T_t}\cA^*\cA\cP_{T_t}}^{-1} \m{Z}_0 } \\
		& \quad + \frobn{ [ \cP_{\tilde{T}_t} ( \cI-\tilde{\cA}^*\tilde{\cA} ) \cP_{\tilde{T}_t} ]^K [ \cP_{\tilde{T}_t}\tilde{\cA}^*\tilde{\cA}\cP_{\tilde{T}_t} ]^{-1} \tilde{\m{Z}}_0 }.
	\end{split}
	\label{eq:W-diff-expansion}
}
The first three terms have be bounded in \eqref{eq:PTminusPT}, \eqref{eq:finite-sums}, and \eqref{eq:I3}, respectively. 
For the last term, we have
\equs{
	\begin{split}
		& \frobn{ [ \cP_{\tilde{T}_t}  (\cI-\tilde{\cA}^*\tilde{\cA}) \cP_{\tilde{T}_t} ]^K [ \cP_{\tilde{T}_t}\tilde{\cA}^*\tilde{\cA}\cP_{\tilde{T}_t} ]^{-1} \tilde{\m{Z}}_0 } \\
		& \le \frac{\delta_{2r+1}^K}{1-\delta_{2r+1}} \frobn{\tilde{\m{Z}}_0} 
		\le \frac{\delta_{2r+1}^K}{1-\delta_{2r+1}} \sbra{ \sqrt{2r}\,\twon{\m{Z}_0} + \frobn{ \m{Z}_0-\tilde{\m{Z}}_0 } } 
		\le \frac{ 9 E_t^2}{\sigma_r(\m{X}_\star)},
	\end{split}
}
where the first inequality follows from Lemma \ref{lem:PT_I_AA_PT},  and the last uses \eqref{eq:Z0}, \eqref{eq:Z0-F}, $\delta_{2r+1} \le 1/2$, and $K=\max\{1,\lceil\log_2 r\rceil\}$.

Consequently, we have
\equ{
	\sup_{(\m{x},\m{y}) \in \cN} \frobn{ \m{W}_{t+1}-\tilde{\m{W}}_{t+1} } 
	\le \frac{26 E_t^2}{\sigma_r(\m{X}_\star)}.
	\label{eq:W-diff}
}

Next, we check the validity of Lemma \ref{lem:TrTr}. Using Weyl's inequality, \eqref{eq:X-next}, $E_t \le c_0 \sigma_r (\m{X}_\star)$, and $\sigma_{r+1} (\m{X}_\star)=0$, we have
\equs{
	\begin{split}
		\sigma_r (\m{W}_{t+1}) 
		& \ge \sigma_r (\m{X}_\star) - \twon{ \m{W}_{t+1} - \m{X}_\star }
		\ge ( 1 - 13 c_0^2 ) \sigma_r (\m{X}_\star), \\
		\sigma_{r+1} (\m{W}_{t+1}) 
		& \le \sigma_{r+1} (\m{X}_\star) + \twon{ \m{W}_{t+1} - \m{X}_\star }
		\le 13 c_0^2  \sigma_r (\m{X}_\star).
	\end{split}
}
Recall that $0<c_0 < 10^{-4}$. It follows that
\equ{
	\begin{split}
		\sigma_r (\m{W}_{t+1}) - \sigma_{r+1} (\m{W}_{t+1}) 
		\ge ( 1 - 26 c_0^2 ) \sigma_r (\m{X}_\star) 
		& \ge \frac{ 1 - 26 c_0^2 }{ 13 c_0^2} \sigma_{r+1} (\m{W}_{t+1}) \\
		& \ge 10^6 \sigma_{r+1} (\m{W}_{t+1}). \label{eq:sigma_r_W}
	\end{split}
}
Moreover,
\equs{
	\begin{split}
		\twon{ \m{W}_{t+1}-\tilde{\m{W}}_{t+1} } 
		& \le \frobn{ \m{W}_{t+1}-\tilde{\m{W}}_{t+1} } 
		\le \frac{ 26 E_t^2}{\sigma_r(\m{X}_\star)} 
		\le 26 c_0^2 \sigma_r(\m{X}_\star) \\
		& \le 10^{-6}( \sigma_r (\m{W}_{t+1}) - \sigma_{r+1} (\m{W}_{t+1}) ),
	\end{split}
}
where the second inequality uses \eqref{eq:W-diff}, the third uses $E_t \le c_0 \sigma_r (\m{X}_\star)$, and the last uses \eqref{eq:sigma_r_W} and $0<c_0 < 10^{-4}$.

Therefore, all the conditions of Lemma \ref{lem:TrTr} are satisfied, and hence we obtain
\equ{
	\begin{split}
		\frobn{\m{X}_{t+1}-\tilde{\m{X}}_{t+1}} 
		& = \frobn{ \cT_r(\m{W}_{t+1}) - \cT_r(\tilde{\m{W}}_{t+1}) } \\	
		& \le (12\times 10^{-6} +20) \frobn{ \m{W}_{t+1}-\tilde{\m{W}}_{t+1} } \\
		& \le  \frac{546 E_t^2}{\sigma_r(\m{X}_\star)},
	\end{split}
	\label{eq:X-virtual}
}
where the last inequality uses \eqref{eq:W-diff}.

\subsection{Putting things together}
\noindent

Finally, combining \eqref{eq:X-next} and \eqref{eq:X-virtual} yields
\equs{
	E_{t+1} 
	\le \frac{ 572 E_t^2}{\sigma_r(\m{X}_\star)},
}

It remains to track the failure probability.
The Gaussian rank-$6r$ RIP event and a union bound over the $K$ required applications of Lemma \ref{lem:spectral-control} and Lemma \ref{lem:ZminusZtilde} show that the preceding estimates fail with probability at most
\equs{
	\exp\sbra{-(n_1+n_2)}
	+
	6 K\exp\sbra{-5(n_1+n_2)}
	\le
	7 \exp\sbra{-(n_1+n_2)},
}
where we use $K \le n_1+n_2 \le \exp((n_1+n_2)/2) $.
This completes the proof.

\bibliographystyle{amsplain}
\bibliography{RGN-sensing}

@article{wu2026lifted,
	title={Lifted Restricted Isometry Property and Fast Algorithm for Spectral Compressed Sensing With Near-Optimal Sample Complexity},
	author={Wu, Xunmeng and Cai, Jian-Feng and Xu, Zongben and Yang, Zai},
	journal={arXiv preprint arXiv:2610.18032},
	year={2026}
}

@article{cai2026near,
	title={Near-Optimal Nonconvex Matrix Completion},
	author={Cai, Jian-Feng and Lu, Xiliang and You, Juntao},
	journal={arXiv preprint arXiv:2609.17048},
	year={2026}
}

@article{huang2017intrinsic,
	title={Intrinsic representation of tangent vectors and vector transports on matrix manifolds},
	author={Huang, Wen and Absil, P-A and Gallivan, Kyle A},
	journal={Numerische Mathematik},
	volume={136},
	number={2},
	pages={523--543},
	year={2017},
	publisher={Springer}
}

@article{breiding2018convergence,
	title={Convergence analysis of {Riemannian Gauss-Newton} methods and its connection with the geometric condition number},
	author={Breiding, Paul and Vannieuwenhoven, Nick},
	journal={Applied Mathematics Letters},
	volume={78},
	pages={42--50},
	year={2018},
	publisher={Elsevier}
}

@article{tanner2013normalized,
	title={Normalized iterative hard thresholding for matrix completion},
	author={Tanner, Jared and Wei, Ke},
	journal={SIAM Journal on Scientific Computing},
	volume={35},
	number={5},
	pages={S104--S125},
	year={2013},
	publisher={SIAM}
}

@article{mohan2012iterative,
	title={Iterative reweighted algorithms for matrix rank minimization},
	author={Mohan, Karthik and Fazel, Maryam},
	journal={The Journal of Machine Learning Research},
	volume={13},
	number={1},
	pages={3441--3473},
	year={2012},
}

@article{fornasier2011low,
	title={Low-rank matrix recovery via iteratively reweighted least squares minimization},
	author={Fornasier, Massimo and Rauhut, Holger and Ward, Rachel},
	journal={SIAM Journal on Optimization},
	volume={21},
	number={4},
	pages={1614--1640},
	year={2011},
	publisher={SIAM}
}

@article{cai2010singular,
	title={A singular value thresholding algorithm for matrix completion},
	author={Cai, Jian-Feng and Cand{\`e}s, Emmanuel J and Shen, Zuowei},
	journal={SIAM Journal on Optimization},
	volume={20},
	number={4},
	pages={1956--1982},
	year={2010},
	publisher={SIAM}
}

@article{wu2024multichannel,
	title={Multichannel frequency estimation with constant amplitude via convex structured low-rank approximation},
	author={Wu, Xunmeng and Yang, Zai and Xu, Zongben},
	journal={SIAM Journal on Matrix Analysis and Applications},
	volume={45},
	number={3},
	pages={1643--1668},
	year={2024},
	publisher={SIAM}
}

@article{chi2019nonconvex,
	title={Nonconvex optimization meets low-rank matrix factorization: An overview},
	author={Chi, Yuejie and Lu, Yue M and Chen, Yuxin},
	journal={IEEE Transactions on Signal Processing},
	volume={67},
	number={20},
	pages={5239--5269},
	year={2019},
	publisher={IEEE}
}

@article{recht2010guaranteed,
	title={Guaranteed minimum-rank solutions of linear matrix equations via nuclear norm minimization},
	author={Recht, Benjamin and Fazel, Maryam and Parrilo, Pablo A},
	journal={SIAM Review},
	volume={52},
	number={3},
	pages={471--501},
	year={2010},
	publisher={SIAM}
}

@article{candes2011robust,
	title={Robust principal component analysis?},
	author={Cand{\`e}s, Emmanuel J and Li, Xiaodong and Ma, Yi and Wright, John},
	journal={Journal of the ACM},
	volume={58},
	number={3},
	pages={1--37},
	year={2011},
	publisher={ACM New York, NY, USA}
}

@article{candes2013phaselift,
	title={Phaselift: Exact and stable signal recovery from magnitude measurements via convex programming},
	author={Candes, Emmanuel J and Strohmer, Thomas and Voroninski, Vladislav},
	journal={Communications on Pure and Applied Mathematics},
	volume={66},
	number={8},
	pages={1241--1274},
	year={2013},
	publisher={Wiley Online Library}
}

@article{candes2012exact,
	title={Exact matrix completion via convex optimization},
	author={Candes, Emmanuel and Recht, Benjamin},
	journal={Communications of the ACM},
	volume={55},
	number={6},
	pages={111--119},
	year={2012},
	publisher={ACM New York, NY, USA}
}

@article{hsu2024quantum,
	title={Quantum state tomography via nonconvex {Riemannian} gradient descent},
	author={Hsu, Ming-Chien and Kuo, En-Jui and Yu, Wei-Hsuan and Cai, Jian-Feng and Hsieh, Min-Hsiu},
	journal={Physical Review Letters},
	volume={132},
	number={24},
	pages={240804},
	year={2024},
	publisher={APS}
}

@article{vandereycken2013low,
	title={Low-rank matrix completion by {Riemannian} optimization},
	author={Vandereycken, Bart},
	journal={SIAM Journal on Optimization},
	volume={23},
	number={2},
	pages={1214--1236},
	year={2013},
	publisher={SIAM}
}

@inproceedings{absil2013extrinsic,
	title={An extrinsic look at the {Riemannian Hessian}},
	author={Absil, P-A and Mahony, Robert and Trumpf, Jochen},
	booktitle={International Conference on Geometric Science of Information},
	pages={361--368},
	year={2013},
	organization={Springer}
}

@book{absil2008optimization,
	title={Optimization algorithms on matrix manifolds},
	author={Absil, P-A and Mahony, Robert and Sepulchre, Rodolphe},
	year={2008},
	publisher={Princeton University Press}
}

@book{vershynin2018high,
	title={High-dimensional probability: An introduction with applications in data science},
	author={Vershynin, Roman},
	volume={47},
	year={2018},
	publisher={Cambridge University Press}
}

@article{candes2011tight,
	title={Tight oracle inequalities for low-rank matrix recovery from a minimal number of noisy random measurements},
	author={Candes, Emmanuel J and Plan, Yaniv},
	journal={IEEE Transactions on Information Theory},
	volume={57},
	number={4},
	pages={2342--2359},
	year={2011},
	publisher={IEEE}
}

@article{stoger2024non,
	title={Non-convex matrix sensing: Breaking the quadratic rank barrier in the sample complexity},
	author={St{\"o}ger, Dominik and Zhu, Yizhe},
	journal={arXiv preprint arXiv:2408.13276},
	year={2024}
}

@inproceedings{tu2016low,
	title={Low-rank solutions of linear matrix equations via {Procrustes} flow},
	author={Tu, Stephen and Boczar, Ross and Simchowitz, Max and Soltanolkotabi, Mahdi and Recht, Ben},
	booktitle={International Conference on Machine Learning},
	pages={964--973},
	year={2016},
	organization={PMLR}
}

@article{tong2021accelerating,
	title={Accelerating ill-conditioned low-rank matrix estimation via scaled gradient descent},
	author={Tong, Tian and Ma, Cong and Chi, Yuejie},
	journal={Journal of Machine Learning Research},
	volume={22},
	number={150},
	pages={1--63},
	year={2021}
}

@article{li2026scaled,
	title={Scaled Gradient Descent for Ill-Conditioned Low-Rank Matrix Recovery with Optimal Sampling Complexity},
	author={Li, Zhenxuan and Huang, Meng},
	journal={arXiv preprint arXiv:2604.00060},
	year={2026}
}

@inproceedings{cai2025fast,
	title={Fast non-convex matrix sensing with optimal sample complexity},
	author={Cai, Jian-Feng and Wu, Tong and Xia, Ruizhe},
	booktitle={The 41st Conference on Uncertainty in Artificial Intelligence},
	year={2025}
}

@article{wei2016guarantees,
	title={Guarantees of {Riemannian} optimization for low rank matrix recovery},
	author={Wei, Ke and Cai, Jian-Feng and Chan, Tony F and Leung, Shingyu},
	journal={SIAM Journal on Matrix Analysis and Applications},
	volume={37},
	number={3},
	pages={1198--1222},
	year={2016},
	publisher={SIAM}
}

@article{zilber2022gnmr,
	title={{GNMR: A provable one-line algorithm for low rank matrix recovery}},
	author={Zilber, Pini and Nadler, Boaz},
	journal={SIAM Journal on Mathematics of Data Science},
	volume={4},
	number={2},
	pages={909--934},
	year={2022},
	publisher={SIAM}
}

@article{kummerle2018harmonic,
	title={Harmonic mean iteratively reweighted least squares for low-rank matrix recovery},
	author={K{\"u}mmerle, Christian and Sigl, Juliane},
	journal={Journal of Machine Learning Research},
	volume={19},
	number={47},
	pages={1--49},
	year={2018}
}

@article{bian2024preconditioned,
	title={A preconditioned {Riemannian} gradient descent algorithm for low-rank matrix recovery},
	author={Bian, Fengmiao and Cai, Jian-Feng and Zhang, Rui},
	journal={SIAM Journal on Matrix Analysis and Applications},
	volume={45},
	number={4},
	pages={2075--2103},
	year={2024},
	publisher={SIAM}
}

@article{luo2023low,
	title={Low-rank tensor estimation via {Riemannian Gauss-Newton}: Statistical optimality and second-order convergence},
	author={Luo, Yuetian and Zhang, Anru R},
	journal={Journal of Machine Learning Research},
	volume={24},
	number={381},
	pages={1--48},
	year={2023}
}

@article{luo2024recursive,
	title={Recursive importance sketching for rank constrained least squares: Algorithms and high-order convergence},
	author={Luo, Yuetian and Huang, Wen and Li, Xudong and Zhang, Anru},
	journal={Operations Research},
	volume={72},
	number={1},
	pages={237--256},
	year={2024},
	publisher={Informs}
}

@article{jain2010guaranteed,
	title={Guaranteed rank minimization via singular value projection},
	author={Jain, Prateek and Meka, Raghu and Dhillon, Inderjit},
	journal={Advances in Neural Information Processing Systems},
	volume={23},
	year={2010}
}

@inproceedings{jain2013low,
	title={Low-rank matrix completion using alternating minimization},
	author={Jain, Prateek and Netrapalli, Praneeth and Sanghavi, Sujay},
	booktitle={Proceedings of the forty-fifth annual ACM symposium on Theory of computing},
	pages={665--674},
	year={2013}
}

\end{document}